\documentclass{aims}
\usepackage{amsmath}
  \usepackage{paralist}
  \usepackage{graphics} 
  \usepackage{epsfig} 
\usepackage{graphicx}  \usepackage{epstopdf} 
 \usepackage[colorlinks=true]{hyperref}
\hypersetup{urlcolor=blue, citecolor=red}

\def\currentvolume{X}
 \def\currentissue{X}
  \def\currentyear{200X}
   \def\currentmonth{XX}
    \def\ppages{X--XX}
     \def\DOI{10.3934/xx.xxxxxxx}

\newtheorem{theorem}{Theorem}[section]
\newtheorem{corollary}{Corollary}

\newtheorem{lemma}[theorem]{Lemma}
\newtheorem{proposition}{Proposition}

\theoremstyle{definition}

\newtheorem{remark}{Remark}

\usepackage{amsmath}
\usepackage{amssymb}
\usepackage{hyperref}
\usepackage{xcolor}

\title[Spectral gaps of growth-fragmentation semigroups]{On spectral gaps of growth-fragmentation semigroups in higher moment
spaces}
\author[Mustapha Mokhtar-Kharroubi and Jacek Banasiak]{}


 \email{mmokhtar@univ-fcomte.fr}
 \email{jacek.banasiak@up.ac.za}

\thanks{Both authors were supported by DSI/NRF SARChI Grant 82770. The second author was also supported by the National Science Centre of Poland Grant 2017/25/B/ST1/00051}

\thanks{$^*$Corresponding author: Mustapha Mokhtar-Kharroubi}
\begin{document}
\maketitle

\centerline{\scshape Mustapha Mokhtar-Kharroubi}
\medskip
{\footnotesize
\centerline{Laboratoire de Math\'{e}matiques, CNRS-UMR 6623}
\centerline{Universit\'{e} de Bourgogne Franche-Comt\'{e}}
\centerline{16 Route de Gray, 25030 Besan\c{c}on, France}}

\medskip

\centerline{\scshape Jacek Banasiak}
\medskip
{\footnotesize
 \centerline{Department of Mathematics and Applied Mathematics}
 \centerline{ University of Pretoria, Pretoria, South Africa}
   \centerline{Institute of Mathematics}
\centerline{{\L\'{o}d\'{z} University of Technology, \L\'{o}d\'{z}, Poland}

}

\bigskip

 \centerline{(Communicated by the associate editor name)}

\begin{abstract}
We present a general approach to proving the existence of spectral gaps and
asynchronous exponential growth for growth-fragmentation semigroups in
moment spaces $L^{1}(%
\mathbb{R}
_{+};\ x^{\alpha }dx)$ and $L^{1}(%
\mathbb{R}
_{+};\ \left( 1+x\right) ^{\alpha }dx)$ for unbounded total fragmentation
rates and continuous growth rates $r(.)$\ such that $\int_{0}^{+\infty }%
\frac{1}{r(\tau )}d\tau =+\infty .\ $The analysis is based on weak
compactness tools and Frobenius theory of positive operators and holds
provided that $\alpha >\widehat{\alpha }$ for a suitable threshold $\widehat{%
\alpha }\geq 1$ that depends on the moment space we consider. A systematic
functional analytic construction is provided. Various examples of
fragmentation kernels illustrating the theory are given and an open problem
is mentioned.
\end{abstract}
\tableofcontents

\section{Introduction}
\subsection{Notation and general assumptions}\mbox{}

This paper deals with the existence of \textit{spectral gaps} (see (\ref%
{Spectral gap}) below) for $C_{0}$-semigroups $\left( V(t)\right) _{t\geq 0}$
governing general growth-fragmentation equations
\begin{eqnarray}
&&\frac{\partial }{\partial t}u(x,t)+\frac{\partial }{\partial x}\left[
r(x)u(x,t)\right] +a(x)u(x,t)  \notag \\
&=&\int_{x}^{+\infty }a(y)b(x,y)u(y,t)dy,\ \ (x,t>0)\
\label{Equation generale}
\end{eqnarray}%
in moment spaces \textit{\ }
\begin{equation}
X_{0,\alpha }:=L^{1}(%
\mathbb{R}
_{+};\ \left( 1+x\right) ^{\alpha }dx),\ \ X_{\alpha }:=L^{1}(%
\mathbb{R}
_{+};\ x^{\alpha }dx)\ \ \ (\alpha >0)  \label{Higher moment spaces}
\end{equation}%
with nonnegative total fragmentation rate
\begin{equation*}
a(\cdot)\in L_{loc}^{1}(0,+\infty )
\end{equation*}%
and a measurable fragmentation kernel$\ b(\cdot,\cdot)$ such that $b(x,y)=0$ if $%
x\geq y$,$\ $
\begin{equation}
\int_{0}^{y}xb(x,y)dx=y  \label{conservativity of mass}
\end{equation}
{and
\begin{equation}
\text{the support of }\left( x,+\infty \right) \ni y\rightarrow a(y)b(x,y)%
\text{ is \textit{unbounded}\ }\left( x>0\right) .
\label{One sided positivity strict}
\end{equation}%
This assumption is required for irreducibility of the growth-fragmentation semigroup, which we also prove under an alternative assumption that  $a(y)>0$ {for} $y\in (0,\infty)$ and  \begin{equation}
\text{there\; is}\; p\in [0,1)\; \text{such\; that\; for\; any}\; y>0,   \inf \mathrm{supp}\; b(\cdot,y)\leq py.
\label{5}
\end{equation}
We note that assumptions \eqref{One sided positivity strict} and \ref{5} have different physical meaning. If \eqref{One sided positivity strict} is satisfied, we allow  for particles of some sizes  not to fragment. This must be offset, however, by the requirement that particles of any size can be obtained by fragmentation of arbitrarily large particles. In this way, non-fragmenting sizes always can be jumped over by daughter particles of parents of a bigger size. Assumption \ref{5} has  different interpretation. The fact that $a(y)>0$ means that particles of any size must split and the second part says that the sizes of daughter particles cannot be too close to the parent's size. In physically realistic situations we expect that fragmentation produces at least two daughter particles whose sizes cannot be both close to the parent's size and thus in such a case this assumption is always satisfied. In particular, a common case of the homogeneous fragmentation kernel $b(x,y) = \frac{1}{y}h\left(\frac{x}{y}\right)$ with $h$ of bounded support, is covered by \eqref{5}. }

Our assumptions on the growth rate are%
\begin{equation}
r\in C(0,+\infty ),\ r(x)>0\ \forall x>0  \label{r est continue et positif}
\end{equation}%
and
\begin{equation}
\int_{0}^{1}\frac{1}{r(\tau )}d\tau <+\infty ,\ \ \int_{1}^{\infty }\frac{1}{%
r(\tau )}d\tau =+\infty   \label{Main assumption}
\end{equation}%
or%
\begin{equation}
\int_{0}^{1}\frac{1}{r(\tau )}d\tau =+\infty ,\ \ \int_{1}^{\infty }\frac{1}{%
r(\tau )}d\tau =+\infty .  \label{Main assumption Bis}
\end{equation}%
In the case (\ref{Main assumption}), we complement (\ref{Equation generale})
with the boundary condition%
\begin{equation*}
\lim_{y\rightarrow 0}r(y)u(y,t)=0.
\end{equation*}%
\ The kinetic equation (\ref{Equation generale}) is the linear part of the growth--coagulation--fragmentation equation where, in the full form, the coagulation part is represented by a quadratic integral term, see \cite{Banasiak-Lamb}. Coagulation and fragmentation processes lay at heart of many fundamental phenomena in ecology, human biology, polymer and aerosol sciences,  astrophysics and the powder production industry; see \cite{BLL} for further details and references, and  \cite{Ber} for the probabilistic context.  A  common feature of these processes is that each involves  a population of inanimate or animate agents that are capable of forming larger or smaller aggregates through, respectively,  coalescence or breakup. Coagulation and fragmentation are conservative processes. In many cases, however, they occur alongside other events that result in the growth of the ensemble. For example, in chemical engineering applications we often observe a precipitation of matter from the solute onto the surface of the aggregates. In biological applications, the growth of the aggregates can occur due to births of new individuals with neonates staying inside the parent's aggregate, see e.g.  \cite{Ackleh1} for the application to phytoplankton or \cite{Okubo} in a general context. The interplay of growth and fragmentation plays also an important role in prion proliferation, see e.g. \cite{Webb}.

In (\ref{Equation generale}), the unknown $u(x,t)$
represents the concentration at time $t$ of \textquotedblleft
agregates\textquotedblright\ with mass\ $x>0$ while $b(x,y)$ $\left(
x<y\right) \ $describes the distribution of mass $x$ aggregates, called
daughter aggregates, spawned by the fragmentation of a mass $y$ aggregate. The
local mass conservation in the fragmentation process is expressed by (\ref%
{conservativity of mass}); we say that the fragmentation kernel $b(\cdot,\cdot)$ is
conservative.

In a preliminary step, we provide explicit formulas of the $C_{0}$%
-semigroups $\left( U(t)\right) _{t\geq 0}$ (with generator $T$)$\ $%
governing the transport equation
\begin{equation}
\frac{\partial }{\partial t}u(x,t)+\frac{\partial }{\partial x}\left[
r(x)u(x,t)\right] +a(x)u(x,t)=0  \label{Growth transsport eq}
\end{equation}%
in the functional spaces $X=X_{\alpha }$ or $X_{0,\alpha }$ and discuss the
effect of the conditions on the growth rate $r(\cdot)$ on them. This direct
approach complements a resolvent Hille-Yosida approach \cite{Banasiak-Lamb};
see also \cite{BLL} Chapter 5 and \cite{BG2}. It turns out, at least for
\textit{bounded} total fragmentation rate $a(\cdot),$ that under Assumption (\ref%
{Main assumption}), the problem (\ref{Growth transsport eq}) is \textit{not}
well-posed in $X_{\alpha }$ in the sense of $C_{0}$-semigroups; see Remark %
\ref{Remark obstruction Bis} below.\ (This does not prevent a generation
theory for suitably \textit{singular} functions $a(\cdot)$ but, in this case,
the whole construction of the paper would need new technicalities; to keep
the coherence of the paper, this special case is treated separately \cite{MK
in prepa}.) Hence, in general, under Assumption (\ref{Main assumption}), a $%
C_{0}$-semigroup $\left( U(t)\right) _{t\geq 0}$ governing (\ref{Growth
transsport eq}) is defined \textit{only} in $X_{0,\alpha }$ and we restrict
our spectral gap construction to this space.$\ $The situation is much more
complex under Assumption (\ref{Main assumption Bis}). Indeed, the problem (%
\ref{Growth transsport eq}) is well-posed in both $X_{\alpha }$ and $%
X_{0,\alpha }$ (with suitable assumptions depending on the space) but the
full spectral gap theory is completed \textit{only} in $X_{\alpha }$.
Indeed, in $X_{0,\alpha }$, although all the preliminary results we need can
be proved, two of them are based on assumptions which are \textit{not}
compatible: indeed, by using the confining role of singular absorptions \cite%
{MK2016}, the resolvent compactness of $T$ (which plays a key role in our
construction) follows from the unboundedness of the total fragmentation rate
$a(\cdot)$ at infinity \textit{and} at zero, while the existence of a $C_{0}$%
-semigroup $\left( V(t)\right) _{t\geq 0}$ with generator%
\begin{equation*}
T+B:D(T)\rightarrow X
\end{equation*}%
($B$ is the fragmentation operator (\ref{Fragmentation operator})) depends
on the boundedness of $a(\cdot)$ near zero. Hence, under Assumption (\ref{Main
assumption Bis}), we need to restrict our construction to the space $%
X_{\alpha }.\ $

In summary, \textit{two spectral gap theories }are given in this paper%
\textit{: }one in $X_{\alpha }$\ under Assumption (\ref{Main assumption Bis}%
) and another one in $X_{0,\alpha }$\ under Assumption (\ref{Main assumption}%
)\textit{. }A spectral gap theory in\textit{\ }$X_{\alpha }\ $under
Assumption (\ref{Main assumption}) needs additional technicalities and is
given in \cite{MK in prepa}. Finally, the existence of a spectral gap in $%
X_{0,\alpha }$ under (\ref{Main assumption Bis}) is an \textit{open problem}%
; see Remark \ref{Remark Open problem}. Our main results are given in
Theorem \ref{Thm Main result} and Theorem \ref{Thm Main result Bis} and are
consequences of many preliminary results of independent interest.
Furthermore, \ various examples of fragmentation kernels (homogeneous or
separable) are given to illustrate the relevance of our assumptions.\

The aim of this paper is twofold. The first aim is well-posedness of (\ref%
{Equation generale}) in the sense of $C_{0}$-semigroups. Indeed, because of
the unboundedness of the total fragmentation rate $a(\cdot)$, the fragmentation
operator (\ref{Fragmentation operator}) is \textit{not} a bounded operator
on $X\ $(where $X$ is either $X_{\alpha }$ or $X_{0,\alpha }).$ Under
suitable assumptions, depending on the space we consider, a generation
result is obtained by using a perturbation theorem by W. Desch specific to
positive semigroups in $L^{1}$-spaces (see below) where the perturbation is
given by the fragmentation operator
\begin{equation}
B: D(T)\ni \varphi \mapsto B\varphi \in X, \quad (B\varphi)(x) := \int_{x}^{+\infty }a(y)b(x,y)\varphi
(y)dy, \label{Fragmentation operator}
\end{equation}%
and $T$ is the generator of $\left( U(t)\right) _{t\geq 0}$ in $X$.$\ $The
second aim of this work is proving the existence of a spectral gap of the $%
C_{0}$-semigroup $\left( V(t)\right) _{t\geq 0}$ governing (\ref{Equation
generale}), i.e., showing that
\begin{equation}
r_{ess}(V(t))<r_{\sigma }(V(t))  \label{Spectral gap}
\end{equation}%
($r_{ess}$ and $r_{\sigma }\ $are respectively the essential spectral radius
and the spectral radius). { If, additionally, $\left(
V(t)\right) _{t\geq 0}$ is  irreducible, then, by \cite{Nagel 1986},  Corollary 3.16 of Chapter C-III,  the spectral bound $\lambda$ of its generator is its dominant eigenvalue and a simple pole of the resolvent. Moreover, by \cite{EN}, Proposition  3.4 of Chapter VI, $\lambda$ is a simple eigenvalue, that is, its eigenspace is one-dimensional. Hence $\left( V(t)\right) _{t\geq 0}$ has the
asynchronous exponential growth property,
\begin{equation}
\left\Vert e^{-\lambda t}V(t)-P\right\Vert _{\mathcal{L}(X)}=O(e^{-%
\varepsilon t})  \label{Asynchronous exponential growth}
\end{equation}%
(for some $\varepsilon >0$), where $P$ is a one-dimensional spectral
projection relative to the isolated algebraically simple dominant eigenvalue
$\lambda $ of the generator, see e.g., \cite{EN}, Theorem 3.5 of Chapter VI, defined as $P = \langle \mathbf{e},\cdot\rangle \mathbf{f}$, where
 $\mathbf{f}$ and $\mathbf{e}$ are strictly positive eigenvectors of, respectively, the generator and its dual, and  $\langle \cdot,\cdot\rangle$ is the  pairing between $X$ and its dual. A summary of these results can also be found in \cite{Vo}, Appendix C. }

  The main mathematical ingredients behind the
occurrence of the spectral gap are a \textit{local} weak compactness property
satisfied by general growth-fragmentation equations (due to the
one-dimensionality of the state variable) and the \textit{confining effect}
of singular total fragmentation rates ensuring the compactness of the
resolvent.

There is an large body of literature dealing with the long term dynamics of solutions to (\ref{Equation generale}) and, in particular, with the asynchronous exponential growth property and the existence of spectral gap. The case when the state space is bounded has been well understood since \cite{Diek}, though the cases with unbounded rates can be tricky, \cite{BPR}. When the state space is unbounded, a number of results have been obtained by the powerful General Relative Entropy method introduced in \cite{MMP}. While the method caters for a large class of coefficients in weighted $L^p$ spaces, the exponential rate of convergence has only been established in \cite{PR}, see also \cite{LP}, and extensively studied since then. Due to its physical interpretation it is important to study the problem in $L^1$ spaces. Some results have been established by probabilistic methods, see e.g., \cite{Ber}\cite{BT} but we are focused on operator--theoretic results for which we refer to the recent works \cite{Mischler}\cite{BG1}\cite{BG2}\cite{Canizo
et al 2020}\cite{Banasiak2018}. In
particular, quantitative estimates of the gap are obtained by means of
Harris's theorem, \cite{Canizo et al 2020}, while \cite{Banasiak2018} contains a comprehensive theory for the discrete case written in the spirit of this paper.\ A special mention should be given to \cite{Doumic-Gabriel}, where the Perron eigenvector and eigenvalue were found and analysed for \eqref{Equation generale} with fairly general coefficients. That paper has stimulated an active research along these lines, culminating in recent works \cite{BG1}\cite{BG2}.

{ We
note also that most of the known literature on spectral gaps deals with Assumption (%
\ref{Main assumption}), see however \cite{Canizo et al 2020}.} Our paper is close in spirit to \cite{BG2} even
if our statements are not the same and our constructions are different and more systematic; see below.

We note that
\begin{equation*}
\int_{0}^{+\infty }u(x,t)xdx,\ \int_{0}^{+\infty }u(x,t)dx,
\end{equation*}%
are respectively the total mass and the total number of agregates at time $%
t\geq 0$. The existence of spectral gaps in the natural functional spaces
\begin{equation*}
X_{1}=L^{1}(%
\mathbb{R}
_{+};\ xdx),\ \ X_{0}=L^{1}(%
\mathbb{R}
_{+};\ dx),\ X_{0,1}=L^{1}(%
\mathbb{R}
_{+};\ \left( 1+x\right) dx),\ \
\end{equation*}%
has been dealt with systematically in \cite{MK2020 Gr-Fragmassloss} but at
the expense of suitable additional \textit{mass loss}%
\begin{equation*}
\int_{0}^{y}xb(x,y)dx=\left( 1-\eta (y)\right) y,\ \ \ (0\leq \eta (y)\leq 1)
\end{equation*}%
or \textit{death} assumptions.\ These assumptions seem to be necessary and
play a key role in well-posedness (via W. Desch's perturbation theorem) of
growth-fragmentation equations in these natural spaces. Fortunately, W.
Desch's theorem can be applied in higher\textit{\ }moment spaces without
such additional assumptions (see \cite{Banasiak-Lamb}, Theorem 2.2).
Actually, we adapt the argument of the proof of (\cite{Banasiak-Lamb},
Theorem 2.2) in our construction. This allows for a significant extension of the
general theory of \cite{MK2020 Gr-Fragmassloss} to higher moment spaces (\ref%
{Higher moment spaces}) by following a similar construction (without
resorting to mass loss or death assumptions) provided that $\alpha >\widehat{%
\alpha }$ for a suitable \textit{threshold}
\begin{equation*}
\widehat{\alpha }\geq 1,
\end{equation*}%
depending on the functional space we consider. This is consistent with the
existence of thresholds known in the literature \cite{Mischler}\cite{BG2}%
\cite{Canizo et al 2020}. As in \cite{MK2020 Gr-Fragmassloss}, our analysis
is based upon few structural assumptions and provides a systematic
functional analytic construction relying on weak compactness tools and the
Frobenius theory of positive operators.

We recall a fundamental perturbation theorem in $L^{1}$ spaces, \cite%
{Desch1988} (see also \cite{Voigt1988}, \cite{MKLivre} Chapter 8 or \cite%
{Banasiak-Arlotti} Chapter 5).

\begin{theorem}
\label{Theorem Desch}(W. Desch's theorem) Let $\left( U(t)\right) _{t\geq 0}$
be a positive $C_{0}$-semigroup on some $L^{1}(\mu )$ space with generator $%
T $ and let $B:D(T)\rightarrow L^{1}(\mu )$ be positive (i.e. $B\varphi \in
L_{+}^{1}(\mu )$ if $\varphi \in L_{+}^{1}(\mu )\cap D(T)$).\ Then
\begin{equation*}
T+B:D(T)\rightarrow L^{1}(\mu )
\end{equation*}%
is a generator of a positive $C_{0}$-semigroup on $L^{1}(\mu )$ if and only
if $T+B$ is resolvent positive or, equivalently, if $\lim_{\lambda
\rightarrow +\infty }r_{\sigma }\left( B(\lambda -T)^{-1}\right) <1.$
\end{theorem}
\subsection{Main results}
\subsubsection{Fully singular growth rates \eqref{Main assumption Bis}}\mbox{}

Let us describe first our main results in the spaces $X_{\alpha }$ and $%
X_{0,\alpha }$ under Assumption (\ref{Main assumption Bis}).\smallskip\\
\textit{Properties of the growth--absorption semigroup in $X_\alpha$.}

The transport $%
C_{0}$-semigroup $\left( U(t)\right) _{t\geq 0}$ governing (\ref{Growth
transsport eq}) exists in $X_{\alpha }$ (resp. in $X_{0,\alpha }$) and is
given by
\begin{equation*}
U(t)f=e^{-\int_{X(y,t)}^{y}\frac{a(p)}{r(p)}dp}f(X(y,t))\frac{\partial X(y,t)%
}{\partial y}
\end{equation*}%
($X(y,t)$ is defined by $\int_{X(y,t)}^{y}\frac{1}{r(\tau )}d\tau =t$)
provided that
\begin{equation}
\varpi :=\sup_{z>0}\frac{r(z)}{z}<+\infty \ \left( \text{resp.\ }\sup_{z>1}\frac{%
r(z)}{z}<+\infty \right).  \label{Growth assumption 1}
\end{equation}%
In addition, Assumptions (\ref{Growth assumption 1}) turn out to be also necessary. Note that under (\ref{Main assumption Bis}), the generation
theory in $X_{0,\alpha }$ needs no condition on the growth rate at the
origin. The resolvent of $T$ is given by
\begin{equation*}
\left( (\lambda -T)^{-1}f\right) (y)=\frac{1}{r(y)}\int_{0}^{y}e^{-%
\int_{x}^{y}\frac{\lambda +a(\tau )}{r(\tau )}d\tau }f(x)dx\ (\Re
\lambda >s(T))
\end{equation*}%
in both spaces $X_{\alpha }$ and $X_{0,\alpha }$.\ \ We show the pointwise
estimate in $X_{\alpha },$
\begin{equation*}
\left\vert (\lambda -T)^{-1}f\right\vert (y)\leq \frac{1}{y^{\alpha }r(y)}%
\left\Vert f\right\Vert _{X_{\alpha }}\ \ (\lambda >\alpha\varpi ).
\end{equation*}%
If we replace the natural condition $\sup_{z>1}\frac{r(z)}{z}<+\infty $ \ by
the stronger one,
\begin{equation}
\widetilde C:=\sup_{z>0}\frac{r(z)}{1+z}<+\infty ,  \label{Stronger assumption}
\end{equation}%
then we can show the pointwise estimate in $X_{0,\alpha },$
\begin{equation}
\left\vert (\lambda -T)^{-1}f\right\vert (y)\leq \frac{1}{\left( 1+y\right)
^{\alpha }r(y)}\left\Vert f\right\Vert _{X_{0,\alpha }}\ \ (\lambda >\alpha \widetilde C).
\label{Poitwise estimate in X zero alpha}
\end{equation}%
We show that $(\lambda - T)^{-1}$ has a smoothing effect in $X_{\alpha }$ in the sense of improving the integrability of the input, that is, for $\lambda
>\alpha\varpi $,%
\begin{equation*}
\int_{0}^{+\infty }\left\vert (\lambda -T)^{-1}f\right\vert (y)a(y)y^{\alpha
}dy\leq \int_{0}^{+\infty }\left\vert f(y)\right\vert y^{\alpha }dy.\
\end{equation*}%
In $X_{0,\alpha }$, \ if we replace the natural condition $\sup_{z>1}\frac{%
r(z)}{z}<+\infty $ by (\ref{Stronger assumption}), we show the smoothing
effect in $X_{0,\alpha }$: for $\lambda >\alpha \widetilde C$%
\begin{equation}
\int_{0}^{+\infty }\left\vert (\lambda -T)^{-1}f\right\vert
(y)a(y)(1+y)^{\alpha }dy\leq \int_{0}^{+\infty }\left\vert f(y)\right\vert
(1+y)^{\alpha }dy.  \label{Smoothing effect iin X zero alpha}
\end{equation}%
The above estimates, combined with the general theory, \cite{MK2016}, on
compactness properties in $L^{1}$ spaces induced by the confining effect of
\textit{singular absorptions}, show that if the sublevel sets of the total
fragmentation rate%
\begin{equation*}
\Omega _{c}=\left\{ x>0;a(x)<c\right\} \ \ \ (c>0)
\end{equation*}%
are ``thin near zero \textit{and} near infinity relatively to $r$" in the
sense%
\begin{equation}
\int_{0}^{+\infty }\frac{1_{\Omega _{c}}(\tau )}{r(\tau )}d\tau <+\infty \ \
(c>0),  \label{Thin sublevel sets}
\end{equation}%
where $1_{\Omega _{c}}$ is the indicator function of $\Omega _{c}\ $(note
that $\frac{1}{r(\cdot)}$ $\notin L^{1}(0,+\infty )$), then $T$ has a compact
resolvent in both spaces $X_{\alpha }$ and $X_{0,\alpha }$. Note that (\ref%
{Thin sublevel sets}) \textit{precludes} $a(\cdot)$ to be bounded near zero or
at infinity. Note also that (\ref{Thin sublevel sets}) occurs for instance
if \ \ \ \
\begin{equation*}
\lim_{y\rightarrow 0^+}a(y)=+\infty ,\ \lim_{y\rightarrow +\infty
}a(y)=+\infty .
\end{equation*}%
\textit{Properties of the full growth--fragmentation semigroup in $X_\alpha$.}

In $X_{\alpha }$, we introduce
\begin{equation*}
n_{\alpha }(y):=\int_{0}^{y}x^{\alpha }b(x,y)dx.
\end{equation*}%
We note that $n_{0}$, abbreviated as
\begin{equation}
n(y):=\int_{0}^{y}b(x,y)dx,  \label{Mean number of daughter}
\end{equation}%
is the mean number (which can be infinite) of daughter aggregates spawned by
the fragmentation of a mass $y$ aggregate.

We show that if
\begin{equation*}
\sup_{y>0}\frac{n_{\alpha }(y)}{y^{\alpha }}<+\infty
\end{equation*}%
(note that it is automatically satisfied if $\alpha\geq 1$), then the fragmentation operator (\ref{Fragmentation operator}) is $T$%
-bounded in $X_{\alpha }$ and
\begin{equation*}
\lim_{\lambda \rightarrow +\infty }\left\Vert B(\lambda -T)^{-1}\right\Vert
_{\mathcal{L}(X_{\alpha })}\leq \lim \sup_{a(y)\rightarrow +\infty }\frac{%
n_{\alpha }(y)}{y^{\alpha }},
\end{equation*}%
where%
\begin{equation*}
\lim \sup_{a(y)\rightarrow +\infty }\frac{n_{\alpha }(y)}{y^{\alpha }}%
:=\lim_{c\rightarrow +\infty }\sup_{\left\{ y;\ a(y)\geq c\right\} }\frac{%
n_{\alpha }(y)}{y^{\alpha }}.
\end{equation*}%
In particular, by W. Desch's perturbation theorem (Theorem \ref{Theorem
Desch}),
\begin{equation*}
T+B:D(T)\subset X_{\alpha }\rightarrow X_{\alpha }\
\end{equation*}%
generates a positive $C_{0}$-semigroup $\left( V(t)\right) _{t\geq 0}$ on $%
X_{\alpha }$ provided that
\begin{equation}
\lim \sup_{a(y)\rightarrow +\infty }\frac{n_{\alpha }(y)}{y^{\alpha }}<1.
\label{Desch condition in X alpha}
\end{equation}
In addition
\begin{equation*}
T+B:D(T)\subset X_{\alpha }\rightarrow X_{\alpha }
\end{equation*}%
is resolvent compact if $T$ is. By exploiting strict comparison results of
spectral radii of positive operators in domination contexts \cite{Marek} and
the convex (weak) compactness property of the strong operator topology \cite%
{Schluchtermann}\cite{MK2004} (see below), we deduce that $\left(
V(t)\right) _{t\geq 0}$ has a spectral gap (\ref{Spectral gap}) and exhibits
the asynchronous exponential growth (\ref{Asynchronous exponential growth})
in $X_{\alpha }$ provided (\ref{One sided positivity strict}) is satisfied
(see Theorem \ref{Thm Main result}). We conjecture that no spectral gap can
occur in $X_{\alpha }$ if $a(\cdot)$ is bounded near zero as suggested by \cite[Theorem 4.1]{BG2}.

One shows the $\alpha $-monotony
\begin{equation*}
\frac{n_{\alpha ^{\prime }}(y)}{y^{\alpha ^{\prime }}}\leq \frac{n_{\alpha
}(y)}{y^{\alpha }}\ (y>0)\ \ (\alpha ^{\prime }>\alpha )
\end{equation*}
{ as well as the $\alpha$-convexity while, obviously, $\frac{n_{1}(y)}{y}= 1.$}

It follows that if $a(\cdot)$ is unbounded at infinity then (\ref{Desch
condition in X alpha}) is \textit{never} satisfied for $\alpha \leq 1.\ $%
{Furthermore, if%
\begin{equation}
\lim \sup_{y
\rightarrow +\infty }\frac{%
n_{\alpha }(y)}{y^{\alpha }}<1  \label{Key point 1}
\end{equation}%
for some $\alpha>1,$ then it is satisfied for all $\alpha>1$ so that if (\ref{Desch condition in X alpha}) is satisfied in some $X_\alpha$, it is satisfied in any $X_\alpha, \alpha>1.$}  Of course, our construction is meaningful only
if (\ref{Key point 1}) holds. {This can be directly checked for instance in the} case of \textit{homogeneous}
fragmentation kernels\textbf{\ }%
\begin{equation}
b(x,y)=\frac{1}{y}h\left(\frac{x}{y}\right)\text{ with }\int_{0}^{1}zh(z)dz=1
\label{Homogeneous kernel}
\end{equation}%
for some%
\begin{equation*}
h\in L_{+}^{1}\left( \left( 0,1\right) ;xdx\right) .\
\end{equation*}%
Indeed, the local conservativeness property%
\begin{equation*}
\int_{0}^{y}xb(x,y)dx=\int_{0}^{y}\frac{x}{y}h\left(\frac{x}{y}%
\right)dx=y\int_{0}^{1}zh(z)dz=y
\end{equation*}%
is satisfied and, for all $\alpha >1$,
\begin{equation*}
\frac{\int_{0}^{y}x^{\alpha }b(x,y)dx}{y^{\alpha }}=y^{-1}\int_{0}^{y}\left(
\frac{x}{y}\right) ^{\alpha }h\left(\frac{x}{y}\right)dx=\int_{0}^{1}z^{\alpha
}h(z)dz<\int_{0}^{1}zh(z)dz=1\
\end{equation*}%
so%
\begin{equation*}
\lim \sup_{a(y)\rightarrow +\infty }\frac{n_{\alpha }(y)}{y^{\alpha }}%
=\int_{0}^{1}z^{\alpha }h(z)dz<1\ (\alpha >1).
\end{equation*}%
We can also check (\ref{Key point 1}) for \textit{separable} (conservative)
fragmentation kernels%
\begin{equation}
b(x,y)=\beta (x)y\left( \int_{0}^{y}s\beta (s)ds\right) ^{-1},
\label{Separable kernel Intro}
\end{equation}%
introduced in \cite{Banasiak 2004}, or even by any convex combination of
such kernels, see Section \ref{Section Separable models}.\smallskip\\
\textit{Analysis in $X_{0,\alpha}$ --- open problems.}

The analysis in $X_{0,\alpha }$ follows the same strategy but the
construction \textit{fails} under Assumption (\ref{Main assumption Bis}).
Indeed, we introduce
\begin{equation}
n_{1,\alpha }(y):=\int_{0}^{y}\left( 1+x\right) ^{\alpha }b(x,y)dx
\label{N 1 alpha}
\end{equation}%
and show that if
\begin{equation}
\sup_{y>0}\frac{n_{1,\alpha }(y)}{\left( 1+y\right) ^{\alpha }}<+\infty ,
\label{Assumption on N1 alpha}
\end{equation}%
then the fragmentation operator (\ref{Fragmentation operator}) is $T$%
-bounded in $X_{0,\alpha }$ and
\begin{equation}
\lim_{\lambda \rightarrow +\infty }\left\Vert B(\lambda -T)^{-1}\right\Vert
_{\mathcal{L}(X_{0,\alpha })}\leq \lim \sup_{a(y)\rightarrow +\infty }\frac{%
n_{1,\alpha }(y)}{\left( 1+y\right) ^{\alpha }}
\label{Limit as lambda tends to infinity}
\end{equation}%
so W. Desch's perturbation theorem shows that if
\begin{equation}
\lim \sup_{a(y)\rightarrow +\infty }\frac{n_{1,\alpha }(y)}{\left(
1+y\right) ^{\alpha }}<1,  \label{Desch condition in X zero alpha}
\end{equation}%
then
\begin{equation}
A:=T+B:D(T)\subset X_{0,\alpha }\rightarrow X_{0,\alpha }\mathbf{\ }
\label{A}
\end{equation}%
generates a positive $C_{0}$-semigroup\textbf{\ }$\left( V(t)\right) _{t\geq
0}\ $on $X_{0,\alpha }.$ Unfortunately, if $a(\cdot)$ is \textit{unbounded near
zero},
\begin{equation*}
\lim \sup_{a(y)\rightarrow +\infty }\frac{n_{1,\alpha }(y)}{\left(
1+y\right) ^{\alpha }}\geq \lim \sup_{y\rightarrow 0}\frac{n_{1,\alpha }(y)}{%
\left( 1+y\right) ^{\alpha }}=\lim \sup_{y\rightarrow 0}n_{1,\alpha }(y)
\end{equation*}%
and (\ref{Desch condition in X zero alpha}) is \textit{not} satisfied since
\begin{equation*}
n_{1,\alpha }(y)\geq \int_{0}^{y}b(x,y)dx=\frac{1}{y}\int_{0}^{y}yb(x,y)dx%
\geq \frac{1}{y}\int_{0}^{y}xb(x,y)dx=1.
\end{equation*}%
On the other hand, under (\ref{Main assumption Bis}), the compactness of the
resolvent of $T$ in $X_{0,\alpha }$ (which plays a key role in our
construction) depends on the unboundedness of $a(\cdot)$ near zero. \textit{This
is why, under Assumption (\ref{Main assumption Bis}), the existence of a
spectral gap in }$X_{0,\alpha }$\textit{\ is an open problem and} \textit{%
our construction is restricted to the space }$X_{\alpha }.$
\subsubsection{Partly singular growth rates \eqref{Main assumption}}\mbox{}

Let us describe now our main results under Assumption (\ref{Main assumption}).\smallskip\\
\textit{Negative results in $X_\alpha$.}

We show first, at least for\textit{\ bounded} total fragmentation
kernels, that (\ref{Growth transsport eq}) is not well-posed in $X_{\alpha }$
in the sense of $C_{0}$-semigroups and consequently, we cannot expect a
generation theory in $X_{\alpha }$ for the full problem (\ref{Equation
generale}) and therefore we restrict ourselves to the space $X_{0,\alpha }.\
$The growth $C_{0}$-semigroup $\left( U(t)\right) _{t\geq 0}$ governing (\ref%
{Growth transsport eq}) with boundary condition%
\begin{equation}
\lim_{x\rightarrow 0}r(x)u(x,t)=0  \label{Boundary condition}
\end{equation}%
exists in the space $X_{0,\alpha }$ and is given by
\begin{equation}
U(t)f=\chi _{\left\{ \int_{0}^{y}\frac{1}{r(\tau )}d\tau >t\right\}
}e^{-\int_{X(y,t)}^{y}\frac{a(p)}{r(p)}dp}f(X(y,t))\frac{\partial X(y,t)}{%
\partial y}  \label{Growth semigroup truncated}
\end{equation}%
($X(y,t)$ is defined by $\int_{X(y,t)}^{y}\frac{1}{r(\tau )}d\tau =t$ for $%
\int_{0}^{y}\frac{1}{r(\tau )}d\tau >t$) provided that (\ref{Stronger
assumption}) is satisfied. This sufficient condition is also ``partly
necessary".\smallskip\\
\textit{Properties of the growth--absorbtion semigroup in $X_{0,\alpha}$.}

We show the estimates (\ref{Poitwise estimate in X zero alpha})(\ref%
{Smoothing effect iin X zero alpha}) in $X_{0,\alpha }$. These estimates,
combined with the general theory \cite{MK2016} on compactness properties in $%
L^{1}$ spaces induced by the confining effect of singular absorptions, show
that if the sublevel sets of the total fragmentation rate
\begin{equation*}
\Omega _{c}=\left\{ x>0;a(x)<c\right\} \ \ \ (c>0)
\end{equation*}%
are ``thin near\textit{\ infinity} relatively to $r$" in the sense%
\begin{equation*}
\int_{1}^{+\infty }\frac{1_{\Omega _{c}}(\tau )}{r(\tau )}d\tau <+\infty \ \
(c>0),
\end{equation*}%
where $1_{\Omega _{c}}$ is the indicator function of $\Omega _{c}\ $(note
that $\frac{1}{r(\cdot)}$ $\notin L^{1}(1,+\infty )$), then $T$ has a compact
resolvent in $X_{0,\alpha }$. This occurs for instance if \ \ \ \
\begin{equation*}
\lim_{y\rightarrow +\infty }a(y)=+\infty .
\end{equation*}%
We introduce (\ref{N 1 alpha}) and show that under (\ref{Assumption on N1
alpha}) the fragmentation operator is $T$-bounded in $X_{0,\alpha }$ and (%
\ref{A}) generates a positive $C_{0}$-semigroup\textbf{\ }$\left(
V(t)\right) _{t\geq 0}\ $on $X_{0,\alpha }$ provided that (\ref{Desch
condition in X zero alpha}) is satisfied. \smallskip\\
\textit{Properties of the full growth--fragmentation semigroup in $X_{0,\alpha}$.}

By restricting ourselves to the case where $a(\cdot)$ is unbounded \textit{at
infinity only}, (\ref{Desch condition in X zero alpha}) amounts to
\begin{equation}
\lim \sup_{y\rightarrow +\infty }\frac{n_{1,\alpha }(y)}{\left( 1+y\right)
^{\alpha }}<1.  \label{Desch condition in X zero alpha Bis}
\end{equation}%
Under (\ref{Desch condition in X zero alpha Bis}),\ the generator (\ref{A})
is resolvent compact in $X_{0,\alpha }$ if $T$ is. By arguing as previously,
we show that $\left( V(t)\right) _{t\geq 0}$ has a spectral gap (\ref%
{Spectral gap}) and exhibits the asynchronous exponential growth (\ref%
{Asynchronous exponential growth}) in $X_{0,\alpha }$ provided (\ref{One
sided positivity strict}) is satisfied (see Theorem \ref{Thm Main result
Bis}).

We show the $\alpha $-monotony
\begin{equation*}
\frac{n_{1,\widehat{\alpha }}(y)}{\left( 1+y\right) ^{\widehat{\alpha }}}%
\leq \frac{n_{1,\alpha }(y)}{\left( 1+y\right) ^{\alpha }}\ (y>0)\ \
(0<\alpha <\widehat{\alpha })
\end{equation*}%
{as well as} the $\alpha$-convexity, and%
\begin{equation*}
\lim \sup_{y\rightarrow +\infty }\frac{n_{1,1}(y)}{ 1+y}\geq 1.
\end{equation*}%
This implies that (\ref{Desch condition in X zero alpha Bis}) is \textit{%
never} satisfied if $\alpha \leq 1.$ It follows that if
\begin{equation}
\lim_{\alpha \rightarrow +\infty }\lim \sup_{y\rightarrow +\infty }\frac{%
n_{1,\alpha }(y)}{\left( 1+y\right) ^{\alpha }}<1,  \label{Key point 2}
\end{equation}%
then there exists a unique threshold
\begin{equation*}
\widetilde{\alpha }:=\inf \left\{ \alpha >1;\lim \sup_{y\rightarrow +\infty }%
\frac{n_{1,\alpha }(y)}{\left( 1+y\right) ^{\alpha }}<1\right\} \geq 1
\end{equation*}%
such that (\ref{Desch condition in X zero alpha Bis}) is satisfied if and
only if $\alpha >\widetilde{\alpha }.$ Similarly, our construction is
meaningfull only if (\ref{Key point 2}) holds. To this end, we show that if
the growth of (\ref{Mean number of daughter}) at infinity is at most
polynomial, i.e., if
\begin{equation}
\eta :=\inf \left\{ \alpha >1;\exists c_{\alpha }>0,\
\int_{0}^{y}b(x,y)dx\leq c_{\alpha }(1+y)^{\alpha }\right\} <+\infty ,
\label{Finite eta  Intro}
\end{equation}%
then%
\begin{equation*}
\lim \sup_{y\rightarrow +\infty }\frac{\int_{0}^{y}\left( 1+x\right)
^{\alpha }b(x,y)dx}{\left( 1+y\right) ^{\alpha }}\leq \lim
\sup_{y\rightarrow +\infty }\frac{\int_{0}^{y}x^{\alpha }b(x,y)dx}{y^{\alpha
}}\ \ (\forall \alpha >\eta )
\end{equation*}%
and, {using $\alpha$-convexity, we show that in this case $1\leq \widetilde{\alpha} \leq  \eta$}.
Again, for a particular example of homogeneous fragmentation kernels (\ref{Homogeneous
kernel})%
\begin{equation*}
n(y)=\int_{0}^{y}\frac{1}{y}h\left(\frac{x}{y}\right)dx=\int_{0}^{1}h(z)dz
\end{equation*}%
and consequently the threshold is exactly one
\begin{equation*}
\widetilde{\alpha }=\eta =1,
\end{equation*}%
provided $\int_{0}^{1}h(z)dz<+\infty$ {(note that here $\eta =1$ by definition, as we only consider exponents bigger than 1)}. We can also check (\ref{Key point 2}%
) for separable (conservative) fragmentation kernels (\ref{Separable kernel
Intro}) or by any convex combination of such kernels, see Section \ref%
{Section Separable models}.
{
We can summarize the results in the following table.
\begin{center}
\begin{tabular}{|c|c|c|}
\hline
Space$\diagdown$ growth& Assumption \eqref{Main assumption}& Assumption \eqref{Main assumption Bis}\\\hline
$X_\alpha$ &$\begin{array}{c}\mathrm{No\; generation\; in\; general}\\
\textrm{Under\; additional}\\\textrm{assumptions,\; see \cite{MK in prepa}}\end{array}$ & $\begin{array}{c}\mathrm{Generation}\\\mathrm{AEG}\end{array}$ \\\hline
$X_{0,\alpha}$& $\begin{array}{c}\mathrm{Generation}\\\mathrm{AEG}\end{array}$& $\begin{array}{c}\mathrm{Generation}\\\mathrm{AEG} - \mathrm{open\, problem}\end{array}$\\\hline
\end{tabular}
\end{center}}
%

\section{The method of characteristics}\mbox{}

{ The explicit formulae for solutions to transport equations (\ref{Growth
transsport eq}) can be obtained by the method of characteristics and belong to the mathematical folklore, see e.g., \cite{Banasiak-Lamb, BG2}; a systematic treatment of them can be found \cite{MK2020 Gr-Fragmassloss}. We recall them for the reader's convenience.}

\begin{proposition}
(\cite[ Proposition 44]{MK2020 Gr-Fragmassloss}) Let (\ref{r est continue et
positif}) and (\ref{Main assumption}) be satisfied. Then the partial
differential equation%
\begin{equation*}
\ \frac{\partial }{\partial t}u(x,t)+\frac{\partial }{\partial x}\left[
r(x)u(x,t)\right] =0,\ (x>0,t>0)
\end{equation*}%
with initial and boundary conditions
\begin{equation*}
u(x,0)=f(x),\ \ \lim_{y\rightarrow 0}r(y)u(y,t)=0\ \ (t>0)
\end{equation*}%
has a unique solution given by%
\begin{equation*}
u(y,t)=\left\{
\begin{array}{c}
\frac{r(X(y,t))f(X(y,t))}{r(y)}=f(X(y,t))\frac{\partial X(y,t)}{\partial y}%
\text{ if \ }\int_{0}^{y}\frac{1}{r(\tau )}d\tau >t \\
0\text{ if }\int_{0}^{y}\frac{1}{r(\tau )}d\tau <t%
\end{array}%
\right.
\end{equation*}%
where $X(y,t)$ is defined, for $\int_{0}^{y}\frac{1}{r(\tau )}d\tau >t$, by%
\begin{equation*}
\int_{X(y,t)}^{y}\frac{1}{r(\tau )}d\tau =t,\ \ X(y,t)\in \left( 0,y\right) .
\end{equation*}
\end{proposition}

\begin{proposition}
(\cite[Proposition 2]{MK2020 Gr-Fragmassloss}) Let (\ref{r est continue et
positif}) and (\ref{Main assumption Bis}) be satisfied. Then the partial
differential equation%
\begin{equation*}
\ \frac{\partial }{\partial t}u(x,t)+\frac{\partial }{\partial x}\left[
r(x)u(x,t)\right] =0,\ (x>0,t>0)
\end{equation*}%
with initial condition $u(x,0)=f(x)$ has a unique solution given by%
\begin{equation*}
u(y,t)=\frac{r(X(y,t))f(X(y,t))}{r(y)},
\end{equation*}%
where $X(y,t)$ ($t>0$) is defined by%
\begin{equation*}
\int_{X(y,t)}^{y}\frac{1}{r(\tau )}d\tau =t,\ \ X(y,t)\in \left( 0,y\right) .
\end{equation*}
\end{proposition}

\section{\label{Section first construction}First construction}\mbox{}

The first construction is based on Assumption (\ref{Main assumption}). It is
devoted to asynchronous exponential growth in the space $X_{0,\alpha }$ only
since we cannot expect in general a generation theory in $X_{\alpha }$ under
(\ref{Main assumption}), see Remark \ref{Remark obstruction} below.

\subsection{Generation theory}\mbox{}

Our first result on the generation of transport semigroups in $X_{0,\alpha }$ is:

\begin{theorem}
\label{Thm generation 1}Let $\alpha >0$ and let (\ref{r est continue et
positif}) and (\ref{Main assumption}) be satisfied. Let $y(x,t)$ be defined
by
\begin{equation}
\int_{x}^{y(x,t)}\frac{1}{r(\tau )}d\tau =t.  \label{Definition y(x,t)}
\end{equation}%
Then $\left( U_{0}(t)\right) _{t\geq 0}$ with
\begin{equation*}
\left( U_{0}(t)f\right) (y)=\chi _{\left\{ \int_{0}^{y}\frac{1}{r(\tau )}%
d\tau >t\right\} }\frac{r(X(y,t))f(X(y,t))}{r(y)}
\end{equation*}%
is a $C_{0}$-semigroup on $X_{0,\alpha }$ if and only if%
\begin{equation*}
\sup_{x>0}\frac{1+y(x,t)}{1+x}<+\infty \ \ (t\geq 0)
\end{equation*}%
and
\begin{equation*}
\left[ 0,+\infty \right) \ni t\mapsto \sup_{x>0}\frac{1+y(x,t)}{1+x}
\end{equation*}%
is locally bounded. In this case%
\begin{equation*}
\left\Vert U_{0}(t)\right\Vert _{\mathcal{L}(X_{0,\alpha })}=\sup_{x>0}\frac{%
\left( 1+y(x,t)\right) ^{\alpha }}{\left( 1+x\right) ^{\alpha }}.
\end{equation*}%
This occurs if there exists $C>0$ such that
\begin{equation}
r(z)\leq C(z+1)\ \ (\forall z>0).  \label{Sublinear}
\end{equation}%
In this case, $\left\Vert U_{0}(t)\right\Vert _{\mathcal{L}(X_{0,\alpha
})}\leq e^{\alpha Ct}.$
\end{theorem}

\begin{proof}
Let us check that $U_{0}(t)$ is a bounded operator on $X_{0,\alpha }.\ $Let $%
y_{0}(t)>0$ be defined by%
\begin{equation}
\int_{0}^{y_{0}(t)}\frac{1}{r(\tau )}d\tau =t.  \label{Definition y0(t)}
\end{equation}%
Note that for $\int_{0}^{y}\frac{1}{r(\tau )}d\tau >t$ we have
\begin{equation}
\int_{X(y,t)}^{y}\frac{1}{r(\tau )}d\tau =t,  \label{Fonction implicit}
\end{equation}%
which shows that (for $t>0$ fixed) $X(y,t)$ is strictly increasing in $y$
and tends to $0$ as $y\rightarrow y_{0}(t)$.\ Note that
\begin{equation*}
(y_{0}(t),+\infty )\ni y\mapsto X(y,t)\in (0,+\infty )
\end{equation*}%
is continuous.\ Note also that (for $t>0$ fixed)
\begin{equation*}
U(y,z):=\int_{z}^{y}\frac{1}{r(\tau )}d\tau -t
\end{equation*}%
is of class $C^{1}$ in $(y,z)$ with
\begin{equation*}
\frac{\partial U(y,z)}{\partial z}=-\frac{1}{r(z)}\neq 0,
\end{equation*}%
so that the implicit function theorem shows that $X(y,t)$ is a $C^{1}$
function in $y\in (y_{0}(t),+\infty )$. Thus, differentiating  (\ref%
{Fonction implicit}) with respect to $y$ we obtain%
\begin{equation*}
\frac{1}{r(y)}-\frac{1}{r(X(y,t))}\frac{\partial X(y,t)}{\partial y}=0
\end{equation*}%
so%
\begin{equation*}
\frac{1}{r(y)}=\frac{1}{r(X(y,t))}\frac{\partial X(y,t)}{\partial y}
\end{equation*}%
and%
\begin{equation*}
\left( U_{0}(t)f\right) (y)=f(X(y,t))\frac{\partial X(y,t)}{\partial y};\ \
y\in (y_{0}(t),+\infty ).
\end{equation*}%
We have
\begin{eqnarray*}
\left\Vert U_{0}(t)f\right\Vert _{X_{0,\alpha }} &=&\int_{0}^{+\infty
}\left\vert \left( U_{0}(t)f\right) (y)\right\vert \left( 1+y\right)
^{\alpha }dy \\
&=&\int_{y_{0}(t)}^{+\infty }\left\vert f(X(y,t))\right\vert \frac{\partial
X(y,t)}{\partial y}\left( 1+y\right) ^{\alpha }dy.
\end{eqnarray*}%
The change of variable $x=X(y,t)$ gives $\ $%
\begin{equation*}
\left\Vert U_{0}(t)f\right\Vert _{X_{0,\alpha }}=\int_{0}^{+\infty
}\left\vert f(x)\right\vert \left( 1+y(x,t)\right) ^{\alpha }dx,
\end{equation*}
where $y(x,t)$ is the unique $y>x$ such that $x=X(y,t)$ i.e. (\ref%
{Definition y(x,t)}). Hence%
\begin{equation*}
\left\Vert U_{0}(t)f\right\Vert _{X_{0,\alpha }}=\int_{0}^{+\infty }\frac{%
\left( 1+y(x,t)\right) ^{\alpha }}{\left( 1+x\right) ^{\alpha }}\left\vert
f(x)\right\vert \left( 1+x\right) ^{\alpha }dx
\end{equation*}%
and $U_{0}(t)$ is a bounded linear operator in $X_{0,\alpha }$ if and only
if
\begin{equation*}
\sup_{x>0}\frac{\left( 1+y(x,t)\right) ^{\alpha }}{\left( 1+x\right)
^{\alpha }}<+\infty .
\end{equation*}%
In such \ a case,
\begin{equation*}
\left\Vert U_{0}(t)\right\Vert _{\mathcal{L}(X_{0,\alpha })}=\sup_{x>0}\frac{%
\left( 1+y(x,t)\right) ^{\alpha }}{\left( 1+x\right) ^{\alpha }}.
\end{equation*}%
Moreover,
\begin{equation*}
\left[ 0,+\infty \right) \ni t\mapsto U_{0}(t)\in \mathcal{L}%
(X_{0,\alpha })
\end{equation*}%
is locally bounded if and only if
\begin{equation*}
\left[ 0,+\infty \right) \ni t\mapsto \sup_{x>0}\frac{\left(
1+y(x,t)\right) ^{\alpha }}{\left( 1+x\right) ^{\alpha }}
\end{equation*}%
is. {Using the flow property of $(y,t)\mapsto X(y,t),$ which follows since it is the solution of an autonomous differential equation, we can prove that the semigroup property  $U_0(s)U_{0}(t) = U_0(t+s), t,s\geq 0,$ is satisfied. It is also easy to see that $ \left[
0,+\infty \right)\ni t \mapsto U_{0}(t)\in \mathcal{L}(X_{0,\alpha })$ is locally bounded. Then, by \cite[Proposition I.1.3]{EN},  to prove that it is a strongly
continuous semigroup on $X_{0,\alpha }$,  it suffices   to check that}%
\begin{equation*}
U_{0}(t)f\rightarrow f\text{ \ in }X_{0,\alpha }\text{ as }t\rightarrow 0
\end{equation*}%
on a dense subspace of $L^{1}(%
\mathbb{R}
_{+};\ \left( 1+x\right) ^{\alpha }dx)$, e.g. for $f$ continuous with
compact support in $(0,+\infty ).$ Note that for any compact set $\left[
c,c^{-1}\right] $%
\begin{equation*}
\int_{0}^{y}\frac{1}{r(\tau )}d\tau >t
\end{equation*}%
for $t$ small enough uniformly in $y\in \left[ c,c^{-1}\right] $ so
\begin{equation*}
\left( U_{0}(t)f\right) (y)=f(X(y,t))\frac{\partial X(y,t)}{\partial y}\ \
\forall y\in \left[ c,c^{-1}\right]
\end{equation*}%
for $t$ small enough. In particular
\begin{equation*}
\int_{X(y,t)}^{y}\frac{1}{r(\tau )}d\tau =t\ \ \ \forall y\in \left[ c,c^{-1}%
\right]
\end{equation*}%
and%
\begin{equation*}
\left( U_{0}(t)f\right) (y)=\frac{r(X(y,t))f(X(y,t))}{r(y)}\ \ \forall y\in %
\left[ c,c^{-1}\right],
\end{equation*}%
for $t$ small enough.\ We note that $X(y,t)\rightarrow y$ as $t\rightarrow 0$
for any $y>0$ and uniformly in\ $y\in \left[ \frac{c}{2},2c^{-1}\right] $.\
Hence
\begin{equation*}
\left( U_{0}(t)f\right) (y)=\frac{r(X(y,t))f(X(y,t))}{r(y)}\rightarrow f(y)\
\ (t\rightarrow 0)
\end{equation*}%
and, by the dominated convergence theorem, $U_{0}(t)f\rightarrow f$ in $%
X_{0,\alpha }$ as $t\rightarrow 0.$

Let us continue with the prove of the last statement of the theorem. For fixed $x>0$, the differentiation in $t\ $of$\ $ $\int_{x}^{y(x,t)}\frac{1%
}{r(\tau )}d\tau =t$ gives%
\begin{equation*}
\frac{1}{r(y(x,t))}\frac{\partial y(x,t)}{\partial t}=1
\end{equation*}%
i.e.%
\begin{equation}
\frac{\partial y(x,t)}{\partial t}=r(y(x,t))\ \ \forall t>0,\text{ with }%
y(x,0)=x.  \label{Ed diff of y}
\end{equation}%
Hence,
\begin{equation*}
y(x,t)=x+\int_{0}^{t}r(y(x,s))ds\leq x+\int_{0}^{t}C\left( y(x,s)+1\right)
ds,
\end{equation*}%
where (\ref{Sublinear}) is used in the last step, so%
\begin{equation*}
y(x,t)+1\leq x+1+\int_{0}^{t}C\left( y(x,s)+1\right) ds
\end{equation*}%
and Gronwall's lemma gives
\begin{equation*}
y(x,t)+1\leq \left( x+1\right) e^{Ct}.
\end{equation*}%
Finally%
\begin{equation*}
\sup_{x>0}\frac{\left( 1+y(x,t)\right) ^{\alpha }}{\left( 1+x\right)
^{\alpha }}\leq e^{\alpha Ct}
\end{equation*}%
and $\left\Vert U_{0}(t)\right\Vert _{\mathcal{L}(X_{0,\alpha })}\leq
e^{\alpha Ct}.$
\end{proof}

\begin{remark}
\label{Remark obstruction}Note that we cannot expect a generation theory in $%
X_{\alpha }$.\ Indeed, $U_{0}(t)$ is bounded in $X_{\alpha }$ if and only
only if
\begin{equation*}
\sup_{x>0}\frac{y(x,t)}{x}<+\infty ,
\end{equation*}%
while (\ref{Definition y(x,t)}) and (\ref{Definition y0(t)}) show that $%
\lim_{x\rightarrow 0}y(x,t)=y_{0}(t)>0.$
\end{remark}

\begin{remark}
\label{Remark partly necessary}Assumption (\ref{Sublinear}) is partly
necessary for the generation theory in $X_{0,\alpha }$, see Remark \ref%
{Remark on optimality} below.\
\end{remark}
{To find the formula of the resolvent of the generator, we take the Laplace transform of $(U_{0}(t))_{t\geqslant 0}.$ We point out that the Laplace integral with respect to $t$ of a continuous $L_1(\mathbb R_+; (1+x)^\alpha dx)$-valued function $t\mapsto f(\cdot,t)$ in the Bochner sense is a.e. in $x$ equal to the Lebesgue integral with respect to $t$ of $f$ treated as a function of two variables $(x,t)\to f(x,t)$, see \cite[Example 2.23]{Banasiak-Arlotti}. Thus, with some change of
variables, we have}
\begin{theorem}\label{Th7}
Let $\alpha >0$, (\ref{r est continue et positif}),\ (\ref{Main assumption})
and (\ref{Sublinear}) be satisfied.\ Let $T_{0}$ be the generator of $%
(U_{0}(t))_{t\geqslant 0}.$ Then
\begin{equation*}
\left( (\lambda -T_{0})^{-1}f\right) (y)=\frac{1}{r(y)}\int_{0}^{y}e^{-%
\int_{x}^{y}\frac{\lambda }{r(s)}ds}f(x)dx,\ \ (f\in X_{0,\alpha })\ (\Re
\lambda >s(T_{0}))
\end{equation*}%
where $s(T_{0})$ is the spectral bound of $T_{0}.$
\end{theorem}
{We note that the last statement follows due to the positivity of $(U_{0}(t))_{t\geqslant 0},$ see \cite[Theorem 1.4.1]{Nerv}.  }
\subsection{A pointwise estimate}\mbox{}

Hereafter we assume that \eqref{Sublinear} is satisfied.

We give the first key a priori estimate.

\begin{lemma}
\label{Lemma pointwise estimate}Let $\alpha >0$, (\ref{r est continue et
positif}), (\ref{Main assumption}) and (\ref{Sublinear}) be satisfied and $%
\lambda \geq \alpha C.$ Then
\begin{equation*}
\left\vert (\lambda -T_{0})^{-1}f\right\vert (y)\leq \frac{1}{\left(
1+y\right) ^{\alpha }r(y)}\left\Vert f\right\Vert _{X_{0,\alpha }}\ \ (f\in
X_{0,\alpha }).
\end{equation*}
\end{lemma}

\begin{proof}
Note that (\ref{Sublinear}), i.e.
\begin{equation*}
\frac{1}{r(\tau )}\geq \frac{C^{-1}}{\tau +1},
\end{equation*}%
implies%
\begin{equation}
e^{-\lambda \int_{x}^{y}\frac{1}{r(\tau )}d\tau }\leq e^{-\frac{\lambda }{C}%
\int_{x}^{y}\frac{1}{\tau +1}d\tau }=e^{-\frac{\lambda }{C}\ln (\frac{y+1}{%
x+1})}=\left(\frac{x+1}{y+1}\right)^{\frac{\lambda }{C}},  \label{Estimate 1}
\end{equation}%
so%
\begin{eqnarray*}
\left\vert (\lambda -T_{0})^{-1}f(y)\right\vert &\leq &\frac{1}{r(y)}%
\int_{0}^{y}e^{-\lambda \int_{x}^{y}\frac{1}{r(\tau )}d\tau }\left\vert
f(x)\right\vert dx \\
&\leq &\frac{1}{r(y)}\int_{0}^{y}\left(\frac{x+1}{y+1}\right)^{\frac{\lambda }{C}%
}\left\vert f(x)\right\vert dx \\
&=&\frac{1}{r(y)}\int_{0}^{y}\frac{1}{\left( 1+x\right) ^{\alpha }}\left(\frac{x+1%
}{y+1}\right)^{\frac{\lambda }{C}}\left\vert f(x)\right\vert \left( 1+x\right)
^{\alpha }dx \\
&=&\frac{1}{\left( 1+y\right) ^{\alpha }r(y)}\int_{0}^{y}\frac{\left(
1+y\right) ^{\alpha }}{\left( 1+x\right) ^{\alpha }}\left(\frac{x+1}{y+1}\right)^{\frac{%
\lambda }{C}}\left\vert f(x)\right\vert \left( 1+x\right) ^{\alpha }dx \\
&=&\frac{1}{\left( 1+y\right) ^{\alpha }r(y)}\int_{0}^{y}\left(\frac{x+1}{y+1}\right)^{%
\frac{\lambda }{C}-\alpha }\left\vert f(x)\right\vert \left( 1+x\right)
^{\alpha }dx.
\end{eqnarray*}%
Finally%
\begin{equation*}
\left\vert (\lambda -T_{0})^{-1}f(y)\right\vert \leq \frac{1}{\left(
1+y\right) ^{\alpha }r(y)}\int_{0}^{y}\left\vert f(x)\right\vert \left(
1+x\right) ^{\alpha }dx\leq \frac{1}{\left( 1+y\right) ^{\alpha }r(y)}%
\left\Vert f\right\Vert _{X_{0,\alpha }}
\end{equation*}%
because $\frac{\lambda }{C}-\alpha \geq 0$ and $\frac{x+1}{y+1}\leq 1$ for $%
0\leq x\leq y.$
\end{proof}

\subsection{The first perturbed semigroup}\mbox{}

We build now a second explicit perturbed $C_{0}$-semigroup by solving, using the method of
characteristics,
\begin{equation*}
\ \frac{\partial }{\partial t}u(x,t)+\frac{\partial }{\partial x}\left[
r(x)u(x,t)\right] +a(x)u(x,t)=0
\end{equation*}%
with initial and boundary conditions%
\begin{equation*}
u(x,0)=f(x),\ \lim_{x\rightarrow 0}r(x)u(x,t)=0.
\end{equation*}%
The solution is given by%
\begin{equation*}
\chi _{\left\{ \int_{0}^{y}\frac{1}{r(\tau )}d\tau >t\right\}
}e^{-\int_{X(y,t)}^{y}\frac{a(p)}{r(p)}dp}\frac{r(X(y,t))f(X(y,t))}{r(y)}.
\end{equation*}%
This defines a perturbed $C_{0}$-semigroup $\left( U(t)\right) _{t\geq 0}$
on $X_{0,\alpha },$ dominated by $\left( U_{0}(t)\right) _{t\geq 0},$
\begin{equation}
\left( U(t)f\right) (y)=\chi _{\left\{ \int_{0}^{y}\frac{1}{r(\tau )}d\tau
>t\right\} }e^{-\int_{X(y,t)}^{y}\frac{a(p)}{r(p)}dp}\frac{r(X(y,t))f(X(y,t))%
}{r(y)}.  \label{First perturbed semigroup}
\end{equation}

\begin{remark}
\label{Remark obstruction Bis}We have seen in Remark \ref{Remark obstruction}
that for $a(\cdot)=0$ we cannot expect a generation theory in $X_{\alpha }.\ $ Hence, by
the bounded perturbation theory, we cannot expect a generation theory in $%
X_{\alpha }$ for bounded $a(\cdot).\ $But this does not prevent (\ref{First
perturbed semigroup}) from defining a $C_{0}$-semigroup in $X_{\alpha }$ for a
suitably singular $a(\cdot).$ Actually, this is the case if $\ a(\cdot)$ is
sufficiently singular at zero but then the whole construction of the paper
needs additional tehnicalities.\ For the sake of clarity, this special case
is treated separately \cite{MK in prepa}.
\end{remark}

As previously, the Laplace transform of $(U(t))_{t\geqslant 0}$ and some
change of variables give:

\begin{proposition}
Let $\alpha >0$ and let (\ref{r est continue et positif}), (\ref{Main
assumption}) and (\ref{Sublinear}) be satisfied. Then, the resolvent of its generator $T$ is given by
\begin{equation*}
\left( \left( \lambda -T\right) ^{-1}f\right) (y)=\frac{1}{r(y)}%
\int_{0}^{y}e^{-\int_{x}^{y}\frac{\lambda +a(\tau )}{r(\tau )}d\tau
}f(x)dx,\ \ \ (f\in X_{0,\alpha })\ (\Re\lambda >s(T)).
\end{equation*}
\end{proposition}
{ As in Theorem \ref{Th7}, the estimate of the abscissa of convergence of the Laplace integral follows from the positivity of $%
(U(t))_{t\geqslant 0}.$}
\subsection{A smoothing effect of the perturbed resolvent}\mbox{}

The second key a priori estimate is given by:

\begin{lemma}
\label{Lemma effet regularisant}Let $\alpha >0$, (\ref{r est continue et
positif}), (\ref{Main assumption}) and (\ref{Sublinear}) be satisfied and$\
\lambda \geq \alpha C.\ $Then
\begin{equation*}
\int_{0}^{+\infty }\left\vert \left( (\lambda -T)^{-1}f\right)
(y)\right\vert a(y)\left( 1+y\right) ^{\alpha }dy\leq \ \int_{0}^{+\infty
}\left\vert (f(y)\right\vert \left( 1+y\right) ^{\alpha }dy,\ \forall f\in
X_{0,\alpha }.
\end{equation*}
\end{lemma}

\begin{proof}
It suffices to consider nonnegative $f.$ Using (\ref{Estimate 1}), we have%
\allowdisplaybreaks
\begin{eqnarray*}
&&\int_{0}^{+\infty }\left( (\lambda -T)^{-1}f\right) (y)a(y)\left(
1+y\right) ^{\alpha }dy \\
&=&\int_{0}^{+\infty }\frac{a(y)\left( 1+y\right) ^{\alpha }}{r(y)}\left(
\int_{0}^{y}e^{-\lambda \int_{x}^{y}\frac{1}{r(p)}dp}e^{-\int_{x}^{y}\frac{%
a(p)}{r(p)}dp}f(x)dx\right) dy \\
&\leq &\int_{0}^{+\infty }\frac{a(y)\left( 1+y\right) ^{\alpha }}{r(y)}%
\left( \int_{0}^{y}\left(\frac{x+1}{y+1}\right)^{\frac{\lambda }{C}}e^{-\int_{x}^{y}%
\frac{a(p)}{r(p)}dp}f(x)dx\right) dy \\
&=&\int_{0}^{+\infty }\left[ \int_{x}^{+\infty }\left(\frac{x+1}{y+1}\right)^{\frac{%
\lambda }{C}}\frac{a(y)\left( 1+y\right) ^{\alpha }}{r(y)}e^{-\int_{x}^{y}%
\frac{a(p)}{r(p)}dp}dy\right] f(x)dx \\
&=&\!\!\!\int_{0}^{+\infty }\!\!\left[ \int_{x}^{+\infty }\frac{1}{\left( 1+x\right)
^{\alpha }}\left(\frac{x+1}{y+1}\right)^{\frac{\lambda }{C}}\frac{a(y)\left( 1+y\right)
^{\alpha }}{r(y)}e^{-\int_{x}^{y}\frac{a(p)}{r(p)}dp}dy\right]\!\! f(x)\left(
1+x\right) ^{\alpha }\!dx \\
&=&\int_{0}^{+\infty }\left[ \int_{x}^{+\infty }\left(\frac{x+1}{y+1}\right)^{\frac{%
\lambda }{C}-\alpha }\frac{a(y)}{r(y)}e^{-\int_{x}^{y}\frac{a(p)}{r(p)}dp}dy%
\right] f(x)\left( 1+x\right) ^{\alpha }dx \\
&\leq &\int_{0}^{+\infty }\left[ \int_{x}^{+\infty }\frac{a(y)}{r(y)}%
e^{-\int_{x}^{y}\frac{a(p)}{r(p)}dp}dy\right] f(x)\left( 1+x\right) ^{\alpha
}dx
\end{eqnarray*}%
where $\frac{x+1}{y+1}\leq 1$ and $\frac{\lambda }{C}-\alpha \geq 0\ $are
used in the last step. Thus%
\begin{eqnarray*}
&&\int_{0}^{+\infty }\left( (\lambda -T)^{-1}f\right) (y)a(y)\left(
1+y\right) ^{\alpha }dy \\
&\leq &\sup_{x>0}\int_{x}^{+\infty }\frac{a(y)}{r(y)}e^{-\int_{x}^{y}\frac{%
a(p)}{r(p)}dp}dy\left( \int_{0}^{+\infty }f(x)\left( 1+x\right) ^{\alpha
}dx\right) .
\end{eqnarray*}%
Finally, the estimate
\begin{equation*}
\begin{split}
\int_{x}^{+\infty }e^{-\int_{x}^{y}\frac{a(p)}{r(p)}dp}\frac{a(y)}{r(y)}%
dy&=-\int_{x}^{+\infty }\frac{d}{dy}\left( e^{-\int_{x}^{y}\frac{a(p)}{r(p)}%
dp}\right) dy\\&=-\left[ e^{-\int_{x}^{y}\frac{a(p)}{r(p)}dp}\right]
_{y=x}^{y=+\infty }\leq 1
\end{split}
\end{equation*}%
ends the proof.
\end{proof}

\subsection{On the full semigroup}\mbox{}

\textbf{\ }We give now the second perturbed semigroup.\qquad

\begin{theorem}
\label{Thm generation full semigroup}Let $\alpha >0$ and let (\ref{r est
continue et positif}), (\ref{Main assumption}) and (\ref{Sublinear}) be
satisfied. Define%
\begin{equation*}
n_{1,\alpha }(y):=\int_{0}^{y}\left( 1+x\right) ^{\alpha }b(x,y)dx.
\end{equation*}%
If
\begin{equation*}
\sup_{y>0}\frac{n_{1,\alpha }(y)}{\left( 1+y\right) ^{\alpha }}<+\infty ,
\end{equation*}%
then the fragmentation operator $B$ is $T$-bounded in $X_{0,\alpha }$ and
\begin{equation*}
\lim_{\lambda \rightarrow +\infty }\left\Vert B(\lambda -T)^{-1}\right\Vert
_{\mathcal{L}(X_{0,\alpha })}\leq \lim \sup_{a(y)\rightarrow +\infty }\frac{%
n_{1,\alpha }(y)}{\left( 1+y\right) ^{\alpha }}.
\end{equation*}%
In particular, if%
\begin{equation}
\lim \sup_{a(y)\rightarrow +\infty }\frac{n_{1,\alpha }(y)}{\left(
1+y\right) ^{\alpha }}<1,  \label{Limit less than 1}
\end{equation}%
then\textbf{\ }%
\begin{equation*}
A:=T+B:X_{0,\alpha }\supset D(T) \rightarrow X_{0,\alpha }
\end{equation*}%
generates a positive $C_{0}$-semigroup\textbf{\ }$\left( V(t)\right) _{t\geq
0}\ $on $X_{0,\alpha }.$
\end{theorem}

\begin{proof}
We note that for nonnegative $\varphi $
\begin{eqnarray*}
\left\Vert B\varphi \right\Vert _{X_{0,\alpha }} &=&\int_{0}^{+\infty
}\left( \int_{x}^{+\infty }a(y)b(x,y)\varphi (y)dy\right) \left( 1+x\right)
^{\alpha }dx \\
&=&\int_{0}^{+\infty }a(y)\left( \int_{0}^{y}\left( 1+x\right) ^{\alpha
}b(x,y)dx\right) \varphi (y)dy \\
&=&\int_{0}^{+\infty }a(y)n_{1,\alpha }(y)\varphi (y)dy.
\end{eqnarray*}%
$\ $Thus, for nonnegative $f$,
\begin{eqnarray}
&&\left\Vert B(\lambda -T)^{-1}f\right\Vert _{X_{0,\alpha }}  \notag \\
&=&\int_{0}^{+\infty }a(y)n_{1,\alpha }(y)\left( (\lambda -T)^{-1}f\right)
(y)dy  \notag \\
&=&\int_{0}^{+\infty }a(y)\frac{n_{1,\alpha }(y)}{\left( 1+y\right) ^{\alpha
}}\left( (\lambda -T)^{-1}f\right) (y)\left( 1+y\right) ^{\alpha }dy.
\label{To be decomposed}
\end{eqnarray}%
Let%
\begin{equation*}
L:=\lim \sup_{a(y)\rightarrow +\infty }\frac{n_{1,\alpha }(y)}{\left(
1+y\right) ^{\alpha }},
\end{equation*}%
that is, for any $\varepsilon >0$ there exists $c_{\varepsilon }>0$ such that%
\begin{equation*}
a(y)\geq c_{\varepsilon }\Longrightarrow \frac{n_{1,\alpha }(y)}{\left(
1+y\right) ^{\alpha }}\leq L+\varepsilon .
\end{equation*}%
We decompose (\ref{To be decomposed}) into two integrals
\begin{eqnarray*}
&&\int_{0}^{+\infty }a(y)\frac{n_{1,\alpha }(y)}{\left( 1+y\right) ^{\alpha }%
}\left( (\lambda -T)^{-1}f\right) (y)\left( 1+y\right) ^{\alpha }dy \\
&=&\int_{\left\{ a(y)\leq c_{\varepsilon }\right\} }a(y)\frac{n_{1,\alpha
}(y)}{\left( 1+y\right) ^{\alpha }}\left( (\lambda -T)^{-1}f\right)
(y)\left( 1+y\right) ^{\alpha }dy \\
&&+\int_{\left\{ a(y)>c_{\varepsilon }\right\} }a(y)\frac{n_{1,\alpha }(y)}{%
\left( 1+y\right) ^{\alpha }}\left( (\lambda -T)^{-1}f\right) (y)\left(
1+y\right) ^{\alpha }dy \\
&=&I_{1}+I_{2}.
\end{eqnarray*}%
We note that
\begin{equation*}
I_{1}\leq c_{\varepsilon }\left\Vert \frac{n_{1,\alpha }(\cdot)}{\left(
1+y\right) ^{\alpha }}\right\Vert _{L^{\infty }}\left\Vert \left( (\lambda
-T)^{-1}f\right) \right\Vert _{X_{0,\alpha }},
\end{equation*}%
while, using Lemma \ref{Lemma effet regularisant},%
\begin{eqnarray*}
I_{2} &\leq &\left( L+\varepsilon \right) \int_{0}^{+\infty }a(y)\left(
(\lambda -T)^{-1}f\right) (y)\left( 1+y\right) ^{\alpha }dy \\
&\leq &\left( L+\varepsilon \right) \left\Vert f\right\Vert _{X_{0,\alpha }}.
\end{eqnarray*}%
Hence,
\begin{eqnarray*}
\left\Vert B(\lambda -T)^{-1}f\right\Vert _{X_{0,\alpha }} &\leq
&c_{\varepsilon }\left\Vert \frac{n_{1,\alpha }(\cdot)}{\left( 1+y\right)
^{\alpha }}\right\Vert _{L^{\infty }}\left\Vert (\lambda -T)^{-1}\right\Vert
_{\mathcal{L}(X_{0,\alpha })}\left\Vert f\right\Vert _{X_{0,\alpha }} \\
&&+\left( L+\varepsilon \right) \left\Vert f\right\Vert _{X_{0,\alpha }}
\end{eqnarray*}%
and%
\begin{equation*}
\left\Vert B(\lambda -T)^{-1}\right\Vert _{\mathcal{L}(X_{0,\alpha })}\leq
c_{\varepsilon }\left\Vert \frac{n_{1,\alpha }(\cdot)}{\left( 1+y\right)
^{\alpha }}\right\Vert _{L^{\infty }}\left\Vert (\lambda -T)^{-1}\right\Vert
_{\mathcal{L}(X_{0,\alpha })}+\left( L+\varepsilon \right) \ \ (\forall
\varepsilon >0).
\end{equation*}%
Since $T$ is a generator of a semigroup, $\left\Vert (\lambda
-T)^{-1}\right\Vert _{\mathcal{L}(X_{0,\alpha })}\rightarrow 0\ $\ as $%
\lambda \rightarrow +\infty $ and hence
\begin{equation*}
\lim_{\lambda \rightarrow +\infty }\left\Vert B(\lambda -T)^{-1}\right\Vert
_{\mathcal{L}(X_{0,\alpha })}\leq L+\varepsilon \ \ (\forall \varepsilon >0).
\end{equation*}%
Consequently,
\begin{equation*}
\lim_{\lambda \rightarrow +\infty }\left\Vert B(\lambda -T)^{-1}\right\Vert
_{\mathcal{L}(X_{0,\alpha })}\leq L.
\end{equation*}%
Finally, if (\ref{Limit less than 1}) is satisfied, $L<1$ and then the
generation follows from the Desch theorem, Theorem \ref{Theorem Desch}.
\end{proof}

\begin{remark}
If $a(\cdot)$ is unbounded near zero, then
\begin{equation*}
\lim \sup_{a(y)\rightarrow +\infty }\frac{n_{1,\alpha }(y)}{\left(
1+y\right) ^{\alpha }}\geq \lim \sup_{y\rightarrow 0}n_{1,\alpha }(y)\geq
\lim \sup_{y\rightarrow 0}\int_{0}^{y}b(x,y)dx\geq 1
\end{equation*}%
because%
\begin{equation*}
\int_{0}^{y}b(x,y)dx=\frac{1}{y}\int_{0}^{y}yb(x,y)dx\geq \frac{1}{y}%
\int_{0}^{y}xb(x,y)dx=1,
\end{equation*}%
so (\ref{Limit less than 1}) is not satisfied. Hence Theorem \ref{Thm
generation full semigroup} is meaningful if $a(\cdot)$ is unbounded only at
infinity; see below.
\end{remark}

Note the useful observation:{
\begin{proposition} Let $(0,\infty)\ni x \to f(x)\in (0,\infty)$ be a non decreasing function and for some there are $y>0, \alpha_0>0$ such that
$$
\int_0^y f^\alpha(x)b(x,y)dx <+\infty$$
for any $\alpha > \alpha_0$. Then
$$
(\alpha_0,+\infty)\ni \alpha \to \frac{\int_0^y f^\alpha(x)b(x,y)dx}{f^\alpha(y)}
$$
is a non increasing and convex function.
\label{Prop monotony in alpha}
\end{proposition}
\begin{proof}
Since $0\leq \frac{f(x)}{f(y)} \leq 1$ for $x\in (0,y],$  $(0,\alpha) \ni \alpha \to \left(\frac{f(x)}{f(y)}\right)^\alpha$ is non increasing and convex, that is, for $0<\alpha_1 \leq \alpha \leq \alpha_2,$
\begin{align*}
\left(\frac{f(x)}{f(y)}\right)^{\alpha_1}&\geq \left(\frac{f(x)}{f(y)}\right)^{\alpha_2},\\
\left(\frac{f(x)}{f(y)}\right)^{\alpha}&\leq \left(\frac{f(x)}{f(y)}\right)^{\alpha_1} + \frac{\left(\frac{f(x)}{f(y)}\right)^{\alpha_2}-\left(\frac{f(x)}{f(y)}\right)^{\alpha_1}}{\alpha_2-\alpha_1}(\alpha -\alpha_1)
\end{align*}
and the statement follows by multiplying both sides with $b(x,y)\geq 0$ and integrating over $(0,y)$ with respect to $x$.
\end{proof}
}
\begin{remark}
 {Applying Proposition \ref{Prop monotony in alpha} to $f(x) = (1+y)^\alpha$ we see that i}f (\ref{Limit less than 1}) is satisfied for some $\alpha >0,$ then it is
satisfied for all $\widehat{\alpha }>\alpha .$
\end{remark}

\begin{remark}
\label{Remark necessity alpha >1}Note first that our assumption (\ref{Limit
less than 1}) precludes the case $\alpha =1\ $if $a(\cdot)$ is unbounded at
infinity. Indeed
\begin{eqnarray*}
n_{1,1}(y) &:=&\int_{0}^{y}\left( 1+x\right)
b(x,y)dx=\int_{0}^{y}b(x,y)dx+\int_{0}^{y}xb(x,y)dx \\
&=&n_{0}(y)+y
\end{eqnarray*}%
and then $\lim \sup_{a(y)\rightarrow +\infty }\frac{n_{1,\alpha }(y)}{\left(
1+y\right) }\geq 1.\ $It follows from Proposition \ref{Prop monotony in
alpha} that
\begin{equation}
\lim \sup_{a(y)\rightarrow +\infty }\frac{n_{1,\alpha }(y)}{\left(
1+y\right) }\geq 1\ \ (0<\alpha \leq 1).\label{n11}
\end{equation}%
Hence the necessity of higher moments, i.e., $\ \alpha >1.\ $More precisely, if%
\begin{equation*}
\lim_{\alpha \rightarrow +\infty }\lim \sup_{a(y)\rightarrow +\infty }\frac{%
n_{1,\alpha }(y)}{\left( 1+y\right) ^{\alpha }}<1,
\end{equation*}%
then the threshold%
\begin{equation}
\widetilde{\alpha }:=\inf \left\{ \alpha >1;\ \lim \sup_{a(y)\rightarrow
+\infty }\frac{n_{1,\alpha }(y)}{\left( 1+y\right) ^{\alpha }}<1\right\}
\label{threshold}
\end{equation}%
is such that (\ref{Limit less than 1}) holds if and only if $\alpha >%
\widetilde{\alpha }.$ See Proposition \ref{Prop estimate threshold} below
for more information about this threshold.
\end{remark}

\begin{remark}
\label{Remark about Lim less than 1}If $a(\cdot)$ is only unbounded at infinity,
then (\ref{Limit less than 1}) amounts to
\begin{equation}
\lim \sup_{y\rightarrow +\infty }\frac{n_{1,\alpha }(y)}{\left( 1+y\right)
^{\alpha }}<1 \label{n1a}
\end{equation}%
and (\ref{threshold}) is given by%
\begin{equation*}
\widetilde{\alpha }:=\inf \left\{ \alpha >1;\ \lim \sup_{y\rightarrow \text{
}+\infty }\frac{n_{1,\alpha }(y)}{\left( 1+y\right) ^{\alpha }}<1\right\}.
\end{equation*}
\end{remark}

We end this subsection by an upper estimate of the threshold (\ref{threshold}%
).\

\begin{proposition}
\label{Prop estimate threshold}We assume that%
\begin{equation}
\eta :=\inf \left\{ \alpha >1;\exists c_{\alpha }>0,\
\int_{0}^{y}b(x,y)dx\leq c_{\alpha }(1+y)^{\alpha }\ \forall y>0\right\}
<+\infty .  \label{Finite eta}
\end{equation}%
Then
\begin{equation*}
\lim \sup_{y\rightarrow +\infty }\frac{\int_{0}^{y}\left( 1+x\right)
^{\alpha }b(x,y)dx}{\left( 1+y\right) ^{\alpha }}\leq \lim
\sup_{y\rightarrow +\infty }\frac{\int_{0}^{y}x^{\alpha }b(x,y)dx}{y^{\alpha
}}\ \ (\forall \alpha >\eta ).
\end{equation*}%
In particular, if \eqref{Finite eta} is satisfied and $a(\cdot)$ is unbounded only at infinity, then  $\widetilde\alpha \leq \eta,$ provided $\widetilde{\alpha}$ is finite.
\end{proposition}

\begin{proof}
Note first that%
\begin{equation*}
\lim \sup_{y\rightarrow +\infty }\frac{\int_{0}^{y}b(x,y)dx}{\left(
1+y\right) ^{\alpha }}=0\ \ (\forall \alpha >\eta ).
\end{equation*}%
Since
\begin{eqnarray*}
\frac{n_{1,\alpha }(y)}{\left( 1+y\right) ^{\alpha }} &=&\frac{%
\int_{0}^{y}\left( 1+x\right) ^{\alpha }b(x,y)dx}{\left( 1+y\right) ^{\alpha
}} \\
&=&\frac{\int_{0}^{y}\zeta (x)b(x,y)dx}{\left( 1+y\right) ^{\alpha }}+\frac{%
\int_{0}^{y}b(x,y)dx}{\left( 1+y\right) ^{\alpha }}
\end{eqnarray*}%
with $\zeta (x)=\left( 1+x\right) ^{\alpha }-1$,%
\begin{equation*}
\lim \sup_{y\rightarrow +\infty }\frac{\int_{0}^{y}\left( 1+x\right)
^{\alpha }b(x,y)dx}{\left( 1+y\right) ^{\alpha }}\leq \lim
\sup_{y\rightarrow +\infty }\frac{\int_{0}^{y}\zeta (x)b(x,y)dx}{\left(
1+y\right) ^{\alpha }}.
\end{equation*}%
Let $\varepsilon >0$ be arbitrary and $M_{\varepsilon }$ be large enough so
that
\begin{equation*}
\frac{\zeta (x)}{x^{\alpha }}\leq 1+\varepsilon \ \ (x\geq M_{\varepsilon }).
\end{equation*}%
Then for $y>M_{\varepsilon }$%
\begin{eqnarray*}
\frac{\int_{0}^{y}\zeta (x)b(x,y)dx}{\left( 1+y\right) ^{\alpha }} &=&\frac{%
\int_{0}^{M_{\varepsilon }}\zeta (x)b(x,y)dx}{\left( 1+y\right) ^{\alpha }}+%
\frac{\int_{M_{\varepsilon }}^{y}\zeta (x)b(x,y)dx}{\left( 1+y\right)
^{\alpha }} \\
&\leq &\zeta (M_{\varepsilon })\frac{\int_{0}^{M_{\varepsilon }}b(x,y)dx}{%
\left( 1+y\right) ^{\alpha }}+\frac{\int_{M_{\varepsilon }}^{y}\zeta
(x)b(x,y)dx}{\left( 1+y\right) ^{\alpha }} \\
&\leq &\zeta (M_{\varepsilon })\frac{\int_{0}^{y}b(x,y)dx}{\left( 1+y\right)
^{\alpha }}+\left( 1+\varepsilon \right) \frac{\int_{0}^{y}x^{\alpha
}b(x,y)dx}{\left( 1+y\right) ^{\alpha }}.
\end{eqnarray*}%
Hence%
\begin{equation*}
\lim \sup_{y\rightarrow +\infty }\frac{\int_{0}^{y}\zeta (x)b(x,y)dx}{\left(
1+y\right) ^{\alpha }}\leq \left( 1+\varepsilon \right) \lim
\sup_{y\rightarrow +\infty }\frac{\int_{0}^{y}x^{\alpha }b(x,y)dx}{\left(
1+y\right) ^{\alpha }}\ \ (\forall \varepsilon >0)
\end{equation*}%
or, equivalently,%
\begin{equation*}
\lim \sup_{y\rightarrow +\infty }\frac{\int_{0}^{y}\zeta (x)b(x,y)dx}{\left(
1+y\right) ^{\alpha }}\leq \left( 1+\varepsilon \right) \lim
\sup_{y\rightarrow +\infty }\frac{\int_{0}^{y}x^{\alpha }b(x,y)dx}{y^{\alpha
}}\ \ (\forall \varepsilon >0).
\end{equation*}%
Hence%
\begin{equation*}
\lim \sup_{y\rightarrow +\infty }\frac{\int_{0}^{y}\zeta (x)b(x,y)dx}{\left(
1+y\right) ^{\alpha }}\leq \lim \sup_{y\rightarrow +\infty }\frac{%
\int_{0}^{y}x^{\alpha }b(x,y)dx}{y^{\alpha }}.
\end{equation*}
{To prove the last statement, we apply Proposition \ref{Prop monotony in alpha}.  If \eqref{n1a} is not satisfied for any $\alpha$, then $\widetilde\alpha =\infty$. If \eqref{n1a} is satisfied for $\alpha <\eta$, then the statement holds. So, we can assume that \eqref{n1a} is satisfied for some $\alpha_2>\eta$. Let us take arbitrary $\alpha_1>\eta$ and $\alpha \in (\alpha_1,\alpha_2).$ By the convexity, for any $y>0$ we have
\begin{align*}
\frac{n_{1,\alpha}(y)}{(1+y)^\alpha} &\leq \frac{n_{1,\alpha_1}(y)}{(1+y)^{\alpha_1}} + \frac{\frac{n_{1,\alpha_2}(y)}{(1+y)^{\alpha_2}} -\frac{n_{1,\alpha_1}(y)}{(1+y)^{\alpha_1}}}{\alpha_2-\alpha_1}(\alpha-\alpha_1)\\
&= \frac{n_{1,\alpha_1}(y)}{(1+y)^{\alpha_1}}\frac{\alpha_2-\alpha}{\alpha_2-\alpha_1} + \frac{n_{1,\alpha_2}(y)}{(1+y)^{\alpha_2}}\frac{\alpha-\alpha_1}{\alpha_2-\alpha_1}.
\end{align*}
For any $\epsilon_1>0$ there is $y_1$ such that for $y>y_1$
$$
\frac{n_{1,\alpha_2}(y)}{(1+y)^{\alpha_2}}\frac{\alpha-\alpha_1}{\alpha_2-\alpha_1}<(1-\epsilon_1)\frac{\alpha-\alpha_1}{\alpha_2-\alpha_1}.
$$
Next, since $\frac{(1+x)^{\alpha_1}}{1+x^{\alpha_1}} \to 1$ as $x\to \infty$, for any $\epsilon_2>0$ we pick $y_2>y_1$ such that
$\frac{(1+x)^{\alpha_1}}{1+x^{\alpha_1}} \leq 1+ \epsilon_2$ for $x\geq y_2$ and,
since $\alpha_1>\eta$, for large $y$ and some $1\leq \eta <\alpha'<\alpha_1$, using \eqref{Finite eta} and
$$
\int_{y_2}^y x^{\alpha_1}b(x,y)dx \leq y^{\alpha_1-1}\int_{y_2}^y x b(x,y)dx \leq y^{\alpha_1-1}\int_{0}^y x b(x,y)dx\leq y^{\alpha_1},
$$
we have
\begin{align*}
&\frac{\int_{0}^{y} (1+x)^{\alpha_1}b(x,y)dx}{(1+y)^{\alpha_1}} \\
&\leq   \frac{\int_{0}^{y_2} (1+x)^{\alpha_1}b(x,y)dx}{(1+y)^{\alpha_1}} +
(1+\epsilon_2)\frac{\int_{y_2}^{y} b(x,y)dx + \int_{y_2}^y x^{\alpha_1}b(x,y)dx}{(1+y)^{\alpha_1}}\\
&\leq   (1+y_2)^{\alpha_1} \frac{c_{\alpha'}(1+y)^{\alpha'}}{(1+y)^{\alpha_1}} +
(1+\epsilon_2)\frac{c_{\alpha'}(1+y)^{\alpha'} + y^{\alpha_1}}{(1+y)^{\alpha_1}}.
\end{align*}
Thus, for any $\epsilon_3$ we have $y_3>y_2$ such that for all $y>y_3$ we have
$$
\frac{\int_{0}^{y} (1+x)^{\alpha_1}b(x,y)dx}{(1+y)^{\alpha_1}} \leq 1+\epsilon_3
$$
and, for such $y$,
\begin{align*}
\frac{n_{1,\alpha}(y)}{(1+y)^\alpha} &\leq (1+\epsilon_3)\frac{\alpha_2-\alpha}{\alpha_2-\alpha_1} +(1-\epsilon_1)\frac{\alpha-\alpha_1}{\alpha_2-\alpha_1}\\
&= 1 + \epsilon_3\frac{\alpha_2-\alpha}{\alpha_2-\alpha_1} - \epsilon_1\frac{\alpha-\alpha_1}{\alpha_2-\alpha_1}.
\end{align*}
Since $\epsilon_3$ and $\epsilon_1$ are independent, taking $\epsilon_3= \frac{\epsilon_1(\alpha-\alpha_1)}{2(\alpha_2-\alpha)}$ we obtain for the corresponding large $y$,
$$
\frac{n_{1,\alpha}(y)}{(1+y)^\alpha} \leq  1- \frac{\epsilon_1}{2}\frac{\alpha-\alpha_1}{\alpha_2-\alpha_1}
$$
hence, since $\alpha_1>\eta$ is arbitrary, we have $\widetilde \alpha \leq \eta$.
}
\end{proof}

\begin{remark}
Similar estimates appear in \cite[Theorem 2.2]{Ban2020}  and \cite[Theorem 2.2]{Banasiak-Lamb}, where W. Desch's theorem is used with the weight $1+x^{\alpha }$
instead of $\left( 1+x\right) ^{\alpha }.\ $
\end{remark}

\begin{remark}
As noted in Introduction, for \textit{homogeneous} fragmentation kernels%
\textbf{\ }%
\begin{equation*}
b(x,y)=\frac{1}{y}h\left(\frac{x}{y}\right)\text{ with }\int_{0}^{1}zh(z)dz=1,
\end{equation*}%
we have%
\begin{equation*}
\frac{\int_{0}^{y}x^{\alpha }b(x,y)dx}{y^{\alpha }}=\int_{0}^{1}z^{\alpha
}h(z)dz<\int_{0}^{1}zh(z)dz=1\ (\alpha >1)
\end{equation*}%
so,
\begin{equation*}
\lim \sup_{y\rightarrow +\infty }\frac{\int_{0}^{y}x^{\alpha }b(x,y)dx}{%
y^{\alpha }}=\int_{0}^{1}z^{\alpha }h(z)dz<\int_{0}^{1}zh(z)dz=1\ (\alpha >1)
\end{equation*}%
and $\widetilde{\alpha }=1.$ \ See Appendix \ref{Section Separable models}
for more examples.
\end{remark}

\subsection{Compactness results}\mbox{}

We start with

\begin{theorem}
\label{Thm T resolvent compact}Let $\alpha >0$,\ (\ref{r est continue et
positif}), (\ref{Main assumption}) and (\ref{Sublinear}) be satisfied. Let
the sublevel sets\break of $a(\cdot)$ be thin at infinity in the sense that for any $%
c>0$
\begin{equation}
\int_{1}^{+\infty }1_{\left\{ a<c\right\} }\frac{1\ }{r(y)}dy<+\infty \
\label{Thin at infinity}
\end{equation}%
(e.g. let $\lim_{x\rightarrow +\infty }a(x)=+\infty $ ). Then $T$ is
resolvent compact.
\end{theorem}

\begin{proof}
Let $\lambda >\alpha C$ and let $f$ \ be in the unit ball of $X_{0,\alpha }$%
, i.e.
\begin{equation*}
\int_{0}^{+\infty }\left\vert f(x)\right\vert \left( 1+x\right) ^{\alpha
}dx\leq 1.
\end{equation*}%
According to Lemma \ref{Lemma effet regularisant}
\begin{equation*}
\int_{0}^{+\infty }\left\vert \left( (\lambda -T)^{-1}f\right)
(y)\right\vert a(y)\left( 1+y\right) ^{\alpha }dy\leq 1.
\end{equation*}%
Let $c>0$ and $\varepsilon >0\ $be \textit{arbitrary}. We have%
\begin{eqnarray*}
1 &\geq &\int_{\varepsilon ^{-1}}^{+\infty }\left\vert \left( (\lambda
-T)^{-1}f\right) (y)\right\vert a(y)\left( 1+y\right) ^{\alpha }dy \\
&=&\int_{\varepsilon ^{-1}}^{+\infty }1_{\left\{ a<c\right\} }\left\vert
\left( (\lambda -T)^{-1}f\right) (y)\right\vert a(y)\left( 1+y\right)
^{\alpha }dy \\
&&+\int_{\varepsilon ^{-1}}^{+\infty }1_{\left\{ a\geq c\right\} }\left\vert
\left( (\lambda -T)^{-1}f\right) (y)\right\vert a(y)\left( 1+y\right)
^{\alpha }dy \\
&\geq &\int_{\varepsilon ^{-1}}^{+\infty }1_{\left\{ a<c\right\} }\left\vert
\left( (\lambda -T)^{-1}f\right) (y)\right\vert a(y)\left( 1+y\right)
^{\alpha }dy+ \\
&&c\int_{\varepsilon ^{-1}}^{+\infty }1_{\left\{ a\geq c\right\} }\left\vert
\left( (\lambda -T)^{-1}f\right) (y)\right\vert \left( 1+y\right) ^{\alpha
}dy,
\end{eqnarray*}%
so%
\begin{equation*}
\sup_{\left\Vert f\right\Vert _{E}\leq 1}\int_{\varepsilon ^{-1}}^{+\infty
}1_{\left\{ a\geq c\right\} }\left\vert \left( (\lambda -T)^{-1}f\right)
(y)\right\vert \left( 1+y\right) ^{\alpha }dy\leq \frac{1}{c}\ \ (\forall
\varepsilon >0).
\end{equation*}%
On the other hand, according to Lemma \ref{Lemma pointwise estimate},
\begin{equation*}
\left\vert (\lambda -T)^{-1}f\right\vert (y)\leq \left\vert (\lambda
-T_{0})^{-1}f\right\vert (y)\leq \frac{1}{\left( 1+y\right) ^{\alpha }r(y)},
\end{equation*}%
so
\begin{equation*}
\int_{\varepsilon ^{-1}}^{+\infty }1_{\left\{ a<c\right\} }\left\vert \left(
(\lambda -T)^{-1}f\right) (y)\right\vert \left( 1+y\right) ^{\alpha }dy\leq
\int_{\varepsilon ^{-1}}^{+\infty }1_{\left\{ a<c\right\} }\frac{1}{r(y)}dy
\end{equation*}%
and then%
\begin{eqnarray*}
&&\int_{\varepsilon ^{-1}}^{+\infty }\left\vert \left( (\lambda
-T)^{-1}f\right) (y)\right\vert \left( 1+y\right) ^{\alpha }dy \\
&=&\int_{\varepsilon ^{-1}}^{+\infty }1_{\left\{ a<c\right\} }\left\vert
\left( (\lambda -T)^{-1}f\right) (y)\right\vert \left( 1+y\right) ^{\alpha
}dy \\
&&+\int_{\varepsilon ^{-1}}^{+\infty }1_{\left\{ a\geq c\right\} }\left\vert
\left( (\lambda -T)^{-1}f\right) (y)\right\vert \left( 1+y\right) ^{\alpha
}dy \\
&\leq &\int_{\varepsilon ^{-1}}^{+\infty }1_{\left\{ a<c\right\} }\frac{1}{%
r(y)}dy+\frac{1}{c}
\end{eqnarray*}%
can be made arbitrarily small (uniformly in $\left\Vert f\right\Vert
_{X_{0,\alpha }}\leq 1$) by choosing \textit{first} $c$ large enough and
then $\varepsilon $ small enough. \

On the other hand on $\left( 0,\varepsilon ^{-1}\right) $ we have
the uniform domination
\begin{equation*}
\left\vert (\lambda -T)^{-1}f\right\vert (y)\leq \ \frac{1_{\left(
0,\varepsilon ^{-1}\right) }(y)}{\left( 1+y\right) ^{\alpha }r(y)}\ \
(\left\Vert f\right\Vert _{X_{0,\alpha }}\leq 1),
\end{equation*}%
where%
\begin{equation*}
\frac{1_{\left( 0,\varepsilon ^{-1}\right) }(y)}{\left( 1+y\right) ^{\alpha
}r(y)}\in X_{0,\alpha }.
\end{equation*}%
Finally $\left\{ (\lambda -T)^{-1}f;\ \left\Vert f\right\Vert _{X_{0,\alpha
}}\leq 1\right\} $ is as close to the relatively weakly compact set
\begin{equation*}
\left\{ 1_{\left( 0,\varepsilon ^{-1}\right) }(y)(\lambda -T)^{-1}f;\
\left\Vert f\right\Vert _{X_{0,\alpha }}\leq 1\right\}
\end{equation*}%
as we want and consequently it is weakly compact. This shows that $(\lambda
-T)^{-1}$ is weakly compact operator and consequently (see \cite[Lemma 14]{MK2016}) $\left( \lambda -T\right) ^{-1}$ is compact.
\end{proof}

\begin{corollary}
\label{Corollary resolvent compact}Let (\ref{r est continue et positif}), (%
\ref{Main assumption}), (\ref{Sublinear}) and (\ref{Limit less than 1}) be
satisfied and let the sublevel sets of $a(\cdot)$ be thin at infinity in the
sense of (\ref{Thin at infinity}). Then $A:=T+B:D(T)\rightarrow X_{0,\alpha
} $, { where $B$ is defined by \eqref{Fragmentation operator},} is resolvent compact.
\end{corollary}

\begin{proof}
This follows simply from Theorem \ref{Thm T resolvent compact} and the fact
that for $\lambda $ large enough, $\sum_{n=0}^{+\infty }\left[ B(\lambda
-T)^{-1}\right] ^{n}$ is a bounded operator and
\begin{equation*}
\left( \lambda -T-B\right) ^{-1}=(\lambda -T)^{-1}\sum_{n=0}^{+\infty }\left[
B(\lambda -T)^{-1}\right] ^{n}.
\end{equation*}
\end{proof}

\subsection{Spectral gap of $(V(t))_{t\geqslant 0}$ in $X_{0,\protect\alpha} $}\mbox{}

{
We start with an irreducibility result which extends \cite[Theorem 5.2.21]{BLL} and is based on ideas from \cite[Proposition 2]{Banasiak2018}.
\begin{theorem}
Let $\alpha >0$, (\ref{r est continue et positif}), (\ref{Main assumption}), (\ref%
{Sublinear}) and (\ref{Limit less than 1}) and either
\begin{enumerate}
\item for any $x>0,$ $\mathrm{supp}_{0\leq y\leq +\infty}\; a(y)b(x, y)$ is unbounded, or
\item $a(y)>0$ for a.a. $y>0$ and there is $0\leq p <1$ such that for each $y>0,$ $\inf \mathrm{supp}\; b(\cdot
, y) \leq py$
\end{enumerate}
 be satisfied. Then $(\lambda
-T-B)^{-1}$ is positivity improving, i.e.,
\begin{equation*}
(\lambda -T-B)^{-1}g>0\text{ a.e.}
\end{equation*}%
for any nontrivial nonnegative $g\in X_{0,\alpha }$ or, equivalently, the $%
C_{0}$-semigroup $\left( V(t)\right) _{t\geq 0}$ is irreducible in $%
X_{0,\alpha }$.
\label{irred}
\end{theorem}}
{\begin{proof} We know that
\begin{equation*}
\left( (\lambda -T)^{-1}g\right) (y)=\frac{1}{r(y)}\int_{0}^{y}e^{-\lambda
\int_{x}^{y}\frac{1}{r(\tau )}d\tau }e^{-\int_{x}^{y}\frac{a(s)}{r(s)}%
ds}g(x)dx
\end{equation*}
and
\begin{equation}
(\lambda -T-B)^{-1}=(\lambda -T)^{-1}\sum_{n=0}^{+\infty }\left[ B(\lambda
-T)^{-1}\right] ^{n},\label{res9}
\end{equation}
Let $g > 0$ and set $z_g = \sup\{z\;: g(z) = 0\; a.a.\; \text{on}\; [0,z]\}$. If $z_g =0$, then obviously $(\lambda -T)^{-1}g>0.$
and the result is valid. Assume then that $z_g >0$ and observe that
$$
\Psi_0(z) := [(\lambda -T)^{-1}g](z) = \left\{\begin{array}{lcl}\frac{1}{r(z)}\int_{z_g}^{z}e^{-
\int_{x}^{z}\frac{\lambda + a(s)}{r(s )}ds }f(x)dx&\text{for}& z\geq z_g,\\
0&\text{for}& 0\leq z < z_g\end{array}\right.
 $$
and $[(\lambda -T)^{-1}g](z)$ is positive for $x>z_g$.   Then, for $x<z_g$,
\begin{align*}
&[(\lambda -T)^{-1}B (\lambda -T)^{-1} g](x) = [(\lambda -T)^{-1}B \Psi_0](x)\\
&=\frac{1}{r(x)}\int_{0}^{x}e^{-
\int_{y}^{x}\frac{\lambda + a(s)}{r(s )}ds }\left(\int_{z_g}^{\infty}a(z)b(y,z)\Psi_0(z)dz\right) dy.
\end{align*}
We see that if assumption 1. is satisfied, then the inner integrand is positive for any $y>0$ and hence $\Psi_1(x):=[(\lambda -T)^{-1}B (\lambda -T)^{-1} g](x)>0$ for any $x>0$. Otherwise, on using assumption 2, $\Psi_1(x)>0$ for $x>\inf \text{supp}\;b(\cdot,z_g).$  Thus, if $p=0$ (in particular, if $b(x,y)>0$ for all $y>0$ and $0< x <y$), then the result is proved. If $p>0$, then $\Psi_1(x)>0$ at least for
 $x>pz_g$. Next, the third term of \eqref{res9} is given by
$$
\Psi_2(x):=[(\lambda -T)^{-1}(B (\lambda -T)^{-1})^2 f](x) = [(\lambda -T)^{-1}B \Psi_1](x)
$$
and thus, by the same argument, $\Psi_2(x)>0$ for $x> p^2z_g$. Using induction and $p^n\to 0$, we conclude that $[(\lambda -T-B)^{-1}f](x)>0$ almost everywhere.
\end{proof}}{
\begin{corollary}
Assumption 2. of Theorem \ref{irred} are satisfied if either
\begin{description}
\item (i) there is $\delta>0$ such that for any $y>0$ we have $n_0(y) \geq 1+\delta$, or
\item (ii) $b\in L_{\infty, loc}(\mathbb R_+\times \mathbb R_+)$.
\end{description}
\end{corollary}}{
\begin{proof}
Assume  that ($i$) is satisfied. If $ \inf \text{supp}\;b(\cdot,y)=0$, then we are done. Otherwise, let for some $y$, $\inf\text{supp}\; b(\cdot, y) = p'y$  for some   $p'\in (0,1)$; $p'$ can depend on $y$. Then
$$
1+\delta \leq n_0(y) = \int_{p'y}^y b(x,y)dx \leq \frac{1}{p'y}\int_{p'y}^y xb(x,y)dx = \frac{1}{p'},
$$
which implies $p'\leq (1+\delta)^{-1}$. Hence, if we select $p$ (independent of $y$) such that $(1+\delta)^{-1} < p<1$, then for each $y>0$, $ \inf \text{supp}\;b(\cdot,y)\leq (1+\delta)^{-1} y< py$. If, instead of ($i$), assumption ($ii$) is satisfied, then the constant $p$ of the proof of Theorem \ref{irred} may be $y$ dependent and though in each step we can prove that the positivity of $\Psi_n$ on $[z_n,\infty)$ implies the positivity of  $\Psi_{n+1}$ on $[z_{n+1}, \infty)$, where $z_{n+1}=p(z_n)z_n$, $z_1 =z_g$, this  sequence may converge (as a decreasing sequence) to a $z_\infty>0$. Then, however, we would have
$$
\text{supp}\; b(\cdot, z_n)\subset (z_\infty,z_n], \quad n \in \mathbb N;
$$
that is,
$$
z_n = \int_{z_\infty}^{z_n} xb(x,z_n)dx.
$$
This, however, leads to a contradiction, since the left-hand side converges to $z_\infty>0$ and the right-hand side, by ($ii$), to 0.
\end{proof}}

\textbf{\ }We are ready to show the main result of the first construction.

\begin{theorem}
\label{Thm Main result} {Under assumptions of Theorem \ref{irred},} let the sublevel
sets of $a(\cdot)$ be thin at infinity in the sense of (\ref{Thin at infinity}).
Then $\left( V(t)\right) _{t\geq 0}$ has a spectral gap, i.e.,%
\begin{equation*}
r_{ess}(V(t))<r_{\sigma }(V(t)),
\end{equation*}%
and it has the asynchronous exponential growth property (\ref{Asynchronous
exponential growth}) in $X_{0,\alpha }$.
\end{theorem}

\begin{proof}
Let%
\begin{equation*}
k(x,y):=1_{\left\{ x<y\right\} }a(y)b(x,y)
\end{equation*}%
be the kernel of $B.\ $ Let further
\begin{equation*}
\overline{k}(x,y):=k(x,y)\wedge 1
\end{equation*}%
and%
\begin{equation*}
\overline{k}_{c}(x,y):=\overline{k}(x,y)p(x)p(y),
\end{equation*}%
where $p\in C(0,+\infty )$ has a compact support in $(0,+\infty )$ and $%
0\leq p(x)\leq 1.$ Note that
\begin{equation*}
k(x,y)\geqslant \overline{k}_{c}(x,y)
\end{equation*}%
and
\begin{eqnarray*}
k(x,y) &=&\left( k(x,y)-\overline{k}_{c}(x,y)\right) +\overline{k}_{c}(x,y)
\\
&=&\widehat{k}(x,y)+\overline{k}_{c}(x,y),
\end{eqnarray*}%
where
\begin{equation*}
\widehat{k}(x,y):=k(x,y)-\overline{k}_{c}(x,y).
\end{equation*}%
Let $\overline{B}$ be the integral operator with kernel $\overline{k}%
_{c}(x,y)$ and let $\widehat{B}$ be the integral operator with kernel $%
\widehat{k}(x,y).$ Since%
\begin{equation*}
\widehat{k}(x,y)\leq k(x,y),
\end{equation*}%
 for $\lambda $ large enough,
\begin{equation*}
\left\Vert \widehat{B}(\lambda -T)^{-1}\right\Vert _{\mathcal{L}(X_{0,\alpha
})}\leq \left\Vert B(\lambda -T)^{-1}\right\Vert _{\mathcal{L}(X_{0,\alpha
})}<1,
\end{equation*}%
so $T+\widehat{B}:D(T)\rightarrow X_{0,\alpha }$ generates a positive
semigroup $(\widehat{V}(t))_{t\geqslant 0}.$ Note that $(V(t))_{t\geqslant
0} $ is generated by
\begin{equation*}
\left( T+\widehat{B}\right) +\overline{B}
\end{equation*}%
where $\overline{B}$ is a \textit{bounded} operator on $X_{0,\alpha }.$
Actually the kernel of $\overline{B}$ is compactly supported in $(0,+\infty
)\times (0,+\infty )$ and bounded and consequently $\overline{B}$ is a
\textit{weakly compact} operator on $X_{0,\alpha }.\ $On the other hand%
\begin{equation*}
V(t)=\widehat{V}(t)+\int_{0}^{t}\widehat{V}(t-s)\overline{B}\widehat{V}(s)ds
\end{equation*}%
and $\int_{0}^{t}\widehat{V}(t-s)\overline{B}\widehat{V}(s)ds$ is a weakly
compact operator (see \cite{Schluchtermann} or \cite{MK2004}) so that $(%
\widehat{V}(t))_{t\geqslant 0}$ and $(V(t))_{t\geqslant 0}$ have the same
essential spectrum \cite{Kato} and then the same essential radius%
\begin{equation}
r_{ess}(\widehat{V}(t))=r_{ess}(V(t)).  \label{Stability essential radius}
\end{equation}%
On the other hand, $\widehat{V}(t)\leq V(t)$,%
\begin{equation*}
(\lambda -T-\widehat{B})^{-1}\leq (\lambda -T-B)^{-1}
\end{equation*}%
and%
\begin{equation*}
(\lambda -T-\widehat{B})^{-1}\neq (\lambda -T-B)^{-1},
\end{equation*}%
because $\overline{B}\neq 0.\ $Since, by Theorem \ref{irred}, $%
(\lambda -T-B)^{-1}$ is positivity improving (and thus irreducible) and
compact (by Corollary \ref{Corollary resolvent compact}),
\begin{equation*}
r_{\sigma }\left[ (\lambda -T-\widehat{B})^{-1}\right] <r_{\sigma }\left[
(\lambda -T-B)^{-1}\right] ,
\end{equation*}%
see \cite{Marek}. Next,
\begin{equation*}
r_{\sigma }\left[ (\lambda -T-\widehat{B})^{-1}\right] =\frac{1}{\lambda
-s(T+\widehat{B})},\ \ r_{\sigma }\left[ (\lambda -T-B)^{-1}\right] =\frac{1%
}{\lambda -s(T+B)}
\end{equation*}%
(see \cite{Nagel 1986}) implies
\begin{equation}
s(T+\widehat{B})<s(T+B)  \label{strict inequality}
\end{equation}%
and hence, in particular,
\begin{equation*}
s(T+B)>-\infty .
\end{equation*}%
Note that the type of a positive semigroup on $L^{1}$-spaces$\ $coincides
with the spectral bound of its generator, see e.g. \cite{Weis}. We combine
this with (\ref{strict inequality}) to get
\begin{equation*}
r_{ess}(V(t))=r_{ess}(\widehat{V}(t))\leq r_{\sigma }(\widehat{V}(t))=e^{s(T+%
\widehat{B})t}<e^{s(T+B)t}=r_{\sigma }(V(t))
\end{equation*}%
so $r_{ess}(V(t))<r_{\sigma }(V(t)).\ ${ Finally, as explained in Introduction, the irreducibility of $\left(
V(t)\right) _{t\geq 0}$ ensures, by \cite[Corollary 3.16 of Chapter C-III]{Nagel 1986}, that the dominant eigenvalue is a simple pole and, by \cite[Proposition  3.4 of Chapter VI]{EN},  we see that it is simple, that is, its eigenspace is one-dimensional.}
\end{proof}

\section{Second construction}\mbox{}

This construction is based on Assumption (\ref{Main assumption Bis}).

\subsection{Generation results}\mbox{}

We start with the space $X_{\alpha }$.

\begin{theorem}
Let \ $\alpha >0$.\ We assume that (\ref{Main assumption Bis}) is satisfied.
Let $X(y,t)$ ($t>0$) be defined by $\int_{X(y,t)}^{y}\frac{1}{r(\tau )}d\tau
=t.$ Then%
\begin{equation*}
\left( U_{0}(t)f\right) (y):=\frac{r(X(y,t))f(X(y,t))}{r(y)}=f(X(y,t))\frac{%
\partial X(y,t)}{\partial y}
\end{equation*}%
defines a positive $C_{0}$-semigroup $(U_{0}(t))_{t\geqslant 0}$ on $%
X_{\alpha }$ if and only if
\begin{equation*}
\ \sup_{x>0}\frac{y(x,t)}{x}<+\infty \ \ (t\geqslant 0)
\end{equation*}%
and
\begin{equation*}
\left[ 0,+\infty \right) \ni t\rightarrow \sup_{x>0}\frac{y(x,t)}{x}\ \
\text{is locally bounded}
\end{equation*}%
where $y(x,t)>x$ is defined by $\ \int_{x}^{y(x,t)}\frac{1}{r(\tau )}d\tau
=t.$ In this case,
\begin{equation*}
\left\Vert U_{0}(t)\right\Vert _{\mathcal{L}(X_{\alpha })}=\sup_{x>0}\frac{%
y^{\alpha }(x,t)}{x^{\alpha }}.
\end{equation*}%
This occurs if
\begin{equation}
C:=\sup_{z>0} \frac{r(z)}{z}<+\infty,  \label{Assumptions on the function r}
\end{equation}%
in which case $\left\Vert U_{0}(t)\right\Vert _{\mathcal{L}(X_{\alpha
})}\leq e^{\alpha Ct}$.
\end{theorem}

\begin{proof}
We set%
\begin{equation*}
U_{0}(t)f:=\frac{r(X(y,t))f(X(y,t))}{r(y)}
\end{equation*}%
and argue as in the proof of Theorem \ref{Thm generation 1}.\ Let us check
that $U_{0}(t)$ is a bounded operator on $X_{\alpha }$.$\ $Note that
\begin{equation}
\int_{X(y,t)}^{y}\frac{1}{r(\tau )}d\tau =t  \label{X(y,t)}
\end{equation}%
and (\ref{Main assumption Bis}) show that (for $t>0$ fixed) $X(y,t)$ is
strictly increasing in $y$ and
\begin{equation*}
\lim_{y\rightarrow 0}X(y,t)=0,\ \lim_{y\rightarrow +\infty }X(y,t)=+\infty .
\end{equation*}%
Since
\begin{equation*}
\frac{1}{r(y)}=\frac{1}{r(X(y,t))}\frac{\partial X(y,t)}{\partial y},
\end{equation*}%
we have
\begin{equation*}
\left( U_{0}(t)f\right) (y)=f(X(y,t))\frac{\partial X(y,t)}{\partial y},\ \
y\in (0,+\infty )
\end{equation*}%
and
\begin{equation*}
\left\Vert U_{0}(t)f\right\Vert _{X_{\alpha }}=\int_{0}^{+\infty }\left\vert
f(X(y,t))\right\vert \frac{\partial X(y,t)}{\partial y}y^{\alpha }dy.
\end{equation*}%
The change of variable $x=X(y,t)$ yields
\begin{equation*}
\left\Vert U_{0}(t)f\right\Vert _{X_{\alpha }}=\int_{0}^{+\infty }\left\vert
f(x)\right\vert y^{\alpha }(x,t)dx,
\end{equation*}%
where $y(x,t)$ is the unique $y>x$ such that $x=X(y,t)$ i.e., $%
\int_{x}^{y(x,t)}\frac{1}{r(\tau )}d\tau =t.$ Since%
\begin{equation*}
\left\Vert U_{0}(t)f\right\Vert _{X_{\alpha }}=\int_{0}^{+\infty }\frac{%
y^{\alpha }(x,t)}{x^{\alpha }}\left\vert f(x)\right\vert x^{\alpha }dx,
\end{equation*}%
$U_{0}(t)$ is a bounded operator on $X_{\alpha }$ if and only if $\sup_{x>0}%
\frac{y(x,t)}{x}<+\infty .$ In this case
\begin{equation*}
\left\Vert U_{0}(t)\right\Vert _{\mathcal{L}(X_{\alpha })}=\sup_{x>0}\frac{%
y^{\alpha }(x,t)}{x^{\alpha }}
\end{equation*}%
and
\begin{equation*}
\left[ 0,+\infty \right) \ni t\rightarrow U_{0}(t)\in \mathcal{L}(X_{\alpha
})
\end{equation*}%
is locally bounded if and only if
\begin{equation*}
\left[ 0,+\infty \right) \ni t\rightarrow \sup_{x>0}\frac{y(x,t)}{x}
\end{equation*}%
is locally bounded.\ { As in the proof of Theorem \ref{Thm generation 1},} to show that $(U_{0}(t))_{t\geqslant 0}$
is strongly continuous on $X_{\alpha }$ it suffices to check that%
\begin{equation*}
U_{0}(t)f\rightarrow f\text{ \ in }L^{1}(%
\mathbb{R}
_{+};\ x^{\alpha }dx)\text{ as }t\rightarrow 0
\end{equation*}%
on a \textit{dense} subspace of $L^{1}(%
\mathbb{R}
_{+};\ x^{\alpha }dx)$, e.g. for $f$ continuous with compact support in $%
(0,+\infty ).$ Note that (\ref{X(y,t)}) shows that $X(y,t)\rightarrow y$ as $%
t\rightarrow 0$ uniformly in $y$ in compact sets of $(0,+\infty )$.\ By
arguing as in the proof of Theorem \ref{Thm generation 1}, one sees that
\begin{equation*}
U_{0}(t)f=\frac{r(X(y,t))f(X(y,t))}{r(y)}\rightarrow f\ \ (t\rightarrow 0)
\end{equation*}%
in $L^{1}(%
\mathbb{R}
_{+};\ x^{\alpha }dx)\ $by the dominated convergence theorem. Finally (\ref%
{Ed diff of y}) implies
\begin{equation}
y(x,t)=x+\int_{0}^{t}r(y(x,s))ds\leq x+\int_{0}^{t}Cy(x,s)ds
\label{Before Gronwall}
\end{equation}%
so, by Gronwall's lemma, $y(x,t)\leq xe^{Ct}$\ and $\sup_{x>0}\frac{%
y^{\alpha }(x,t)}{x^{\alpha }}\leq e^{\alpha Ct}.$
\end{proof}

\begin{remark}
\label{Remark on optimality}One can show (see \cite[Proposition 6]{MK2020 Gr-Fragmassloss}) that if%
\begin{equation*}
\lim_{z\rightarrow 0}\frac{r(z)}{z}=+\infty \ \text{or}\lim_{z\rightarrow
+\infty }\frac{r(z)}{z}=+\infty,
\end{equation*}%
then $\sup_{x>0}\frac{y(x,t)}{x}=+\infty .\ $In particular we have not a
generation theory in $X_{\alpha }.$ This shows the ``optimality" of
Assumption (\ref{Assumptions on the function r}) in $X_{\alpha }.\ $This
shows also\ that in Theorem \ref{Thm generation 1}, (\ref{Sublinear}) is
partly necessary.
\end{remark}

We deal now with $X_{0,\alpha }.$

\begin{theorem}
Let \ $\alpha >0$.\ We assume that (\ref{Main assumption Bis}) is satisfied.
Let $X(y,t)$ ($t>0$) be defined by $\int_{X(y,t)}^{y}\frac{1}{r(\tau )}d\tau
=t.$ Then%
\begin{equation*}
\left( U_{0}(t)f\right) (y):=\frac{r(X(y,t))f(X(y,t))}{r(y)}=f(X(y,t))\frac{%
\partial X(y,t)}{\partial y}
\end{equation*}%
defines a positive $C_{0}$-semigroup $(U_{0}(t))_{t\geqslant 0}$ on $%
X_{0,\alpha }$ if and only if
\begin{equation*}
\ \sup_{x>0}\frac{1+y(x,t)}{1+x}<+\infty \ \ (t\geqslant 0)
\end{equation*}%
and
\begin{equation*}
\left[ 0,+\infty \right) \ni t\rightarrow \sup_{x>0}\frac{1+y(x,t)}{1+x}\ \
\text{is locally bounded,}
\end{equation*}%
where $y(x,t)\geq x$ is defined by $\int_{x}^{y(x,t)}\frac{1}{r(\tau )}d\tau
=t.$ In this case
\begin{equation*}
\left\Vert U_{0}(t)\right\Vert _{\mathcal{L}(X_{0,\alpha })}=\sup_{x>0}\frac{%
\left( 1+y(x,t)\right) ^{\alpha }}{\left( 1+x\right) ^{\alpha }}\text{.}
\end{equation*}%
This occurs if
\begin{equation}
\widehat{C}:=\sup_{z>1}\frac{r(z)}{z}<+\infty .  \label{Assumption on r bis}
\end{equation}
\end{theorem}

\begin{proof}
Arguing as in the previous proof, we obtain%
\begin{equation*}
\left\Vert U_{0}(t)f\right\Vert _{X_{0,\alpha }}=\int_{0}^{+\infty
}\left\vert f(X(y,t))\right\vert \frac{\partial X(y,t)}{\partial y}\left(
1+y\right) ^{\alpha }dy,
\end{equation*}%
so%
\begin{eqnarray*}
\left\Vert U_{0}(t)f\right\Vert _{X_{0,\alpha }} &=&\int_{0}^{+\infty
}\left\vert f(x)\right\vert \left( 1+y(x,t)\right) ^{\alpha }dx \\
&=&\int_{0}^{+\infty }\frac{\left( 1+y(x,t)\right) ^{\alpha }}{\left(
1+x\right) ^{\alpha }}\left\vert f(x)\right\vert \left( 1+x\right) ^{\alpha
}dx
\end{eqnarray*}%
and%
\begin{equation*}
\left\Vert U_{0}(t)\right\Vert _{\mathcal{L}(X_{0,\alpha })}=\sup_{x>0}\frac{%
\left( 1+y(x,t)\right) ^{\alpha }}{\left( 1+x\right) ^{\alpha }}.
\end{equation*}%
Note that $\int_{x}^{y(x,t)}\frac{1}{r(\tau )}d\tau =t$ implies that $%
\lim_{x\rightarrow 0}y(x,t)=0$ uniformly for bounded sets of $t$, so, for any $%
\overline{t}>0$,
\begin{equation*}
\sup_{t\in \left[ 0,\overline{t}\right] }\sup_{x<1}\frac{1+y(x,t)}{\left(
1+x\right) }<+\infty .
\end{equation*}%
Since
\begin{equation*}
y(x,t)=x+\int_{0}^{t}r(y(x,s))ds,
\end{equation*}%
 $y(x,t)\geq x$ and, by (\ref{Assumption on r bis}),
\begin{equation*}
y(x,t)\leq x+\int_{0}^{t}\widehat{C}y(x,s)ds\ \ (x>1).
\end{equation*}%
Hence,
\begin{equation*}
1+y(x,t)\leq 1+x+\int_{0}^{t}\widehat{C}\left( 1+y(x,s)\right) ds\ (x>1)
\end{equation*}%
and%
\begin{equation*}
1+y(x,t)\leq \left( 1+x\right) e^{\widehat{C}t}\ \ (x>1)
\end{equation*}%
by Gronwall's inequality.\ Finally,
\begin{equation*}
t\rightarrow \sup_{x>0}\frac{\left( 1+y(x,t)\right) ^{\alpha }}{\left(
1+x\right) ^{\alpha }}<+\infty
\end{equation*}%
is locally bounded.\ The rest of the proof is the same as the previous one.
\end{proof}

\begin{remark}
As in Remark \ref{Remark on optimality}, if%
\begin{equation*}
\ \lim_{z\rightarrow +\infty }\frac{r(z)}{z}=+\infty,
\end{equation*}%
then $\sup_{x>0}\frac{y(x,t)}{x}=+\infty .\ $This again shows the ``optimality" of
Assumption (\ref{Assumption on r bis}) in $X_{0,\alpha }$.
\end{remark}

\subsection{A pointwise estimate}\mbox{}

We give now the first a priori estimate in the spaces$\ X_{\alpha }$ and $%
X_{0,\alpha }.$

\begin{lemma}
\label{Lemma pointwise estimate Bis}Let $\alpha >0$\ and let (\ref{Main
assumption Bis}) be satisfied.

(i) Let (\ref{Assumptions on the function r}) be satisfied and $\lambda \geq
\alpha C.$ Then
\begin{equation*}
\left\vert (\lambda -T_{0})^{-1}f\right\vert (y)\leq \frac{1}{y^{\alpha }r(y)%
}\left\Vert f\right\Vert _{X_{\alpha }}\ \ (f\in X_{\alpha }).
\end{equation*}

(ii) Let (\ref{Stronger assumption}) be satisfied and $\lambda \geq \alpha
\widetilde{C}$. Then
\begin{equation*}
\left\vert (\lambda -T_{0})^{-1}f\right\vert (y)\leq \frac{1}{\left(
1+y\right) ^{\alpha }r(y)}\left\Vert f\right\Vert _{X_{0,\alpha }}\ \ (f\in
X_{0,\alpha }).
\end{equation*}
\end{lemma}

\begin{proof}
(i) Note that $r(\tau )\leq C\tau \ \ (\forall \tau >0)$, that is,
\begin{equation*}
\frac{1}{r(\tau )}\geq \frac{1}{C\tau }
\end{equation*}%
implies%
\begin{equation}
e^{-\lambda \int_{x}^{y}\frac{1}{r(\tau )}d\tau }\leq e^{-\frac{\lambda }{C}%
\int_{x}^{y}\frac{1}{\tau }d\tau }=e^{-\frac{\lambda }{C}\ln \left(\frac{y}{x}\right)}=\left(
\frac{x}{y}\right)^{\frac{\lambda }{C}},  \label{Estimate 2}
\end{equation}%
so%
\begin{eqnarray*}
\left\vert (\lambda -T_{0})^{-1}f(y)\right\vert &\leq &\frac{1}{r(y)}%
\int_{0}^{y}e^{-\lambda \int_{x}^{y}\frac{1}{r(\tau )}d\tau }\left\vert
f(x)\right\vert dx \\
&\leq &\frac{1}{r(y)}\int_{0}^{y}\left(\frac{x}{y}\right)^{\frac{\lambda }{C}%
}\left\vert f(x)\right\vert dx.
\end{eqnarray*}%
Since $f\in X_{\alpha }$,%
\begin{eqnarray*}
\left\vert (\lambda -T_{0})^{-1}f(y)\right\vert &\leq &\frac{1}{r(y)}%
\int_{0}^{y}x^{-\alpha }\left(\frac{x}{y}\right)^{\frac{\lambda }{C}}\left\vert
f(x)\right\vert x^{\alpha }dx \\
&=&\frac{1}{y^{\alpha }r(y)}\int_{0}^{y}y^{\alpha }x^{-\alpha }\left(\frac{x}{y}%
\right)^{\frac{\lambda }{C}}\left\vert f(x)\right\vert x^{\alpha }dx \\
&=&\frac{1}{y^{\alpha }r(y)}\int_{0}^{y}\left(\frac{x}{y}\right)^{\frac{\lambda }{C}%
-\alpha }\left\vert f(x)\right\vert x^{\alpha }dx \\
&\leq &\frac{1}{y^{\alpha }r(y)}\int_{0}^{y}\left\vert f(x)\right\vert
x^{\alpha }dx\leq \frac{1}{y^{\alpha }r(y)}\left\Vert f\right\Vert
_{X_{\alpha }},
\end{eqnarray*}%
on account of $\frac{x}{y}\leq 1$ and $\frac{\lambda }{C}-\alpha \geq 0.$

(ii) Note that%
\begin{equation*}
\frac{1}{r(\tau )}\geq \frac{1}{\widetilde{C}(\tau +1)}\ (\tau>0)
\end{equation*}%
and
\begin{equation}
e^{-\lambda \int_{x}^{y}\frac{1}{r(\tau )}d\tau }\leq e^{-\frac{\lambda }{%
\widetilde{C}}\int_{x}^{y}\frac{1}{\tau +1}d\tau }=e^{-\frac{\lambda }{%
\widetilde{C}}\ln \left(\frac{y+1}{x+1}\right)}=\left(\frac{x+1}{y+1}\right)^{\frac{\lambda }{%
\widetilde{C}}}  \label{Estimate 2 Bis}
\end{equation}%
so that if $f\in X_{0,\alpha },$ then
\begin{eqnarray*}
\left\vert (\lambda -T_{0})^{-1}f(y)\right\vert &\leq &\frac{1}{r(y)}%
\int_{0}^{y}e^{-\lambda \int_{x}^{y}\frac{1}{r(\tau )}d\tau }\left\vert
f(x)\right\vert dx \\
&\leq &\frac{1}{r(y)}\int_{0}^{y}\left(\frac{x+1}{y+1}\right)^{\frac{\lambda }{\widetilde{
C}}}\left\vert f(x)\right\vert dx \\
&=&\frac{1}{r(y)}\int_{0}^{y}\frac{1}{\left( 1+x\right) ^{\alpha }}\left(\frac{x+1%
}{y+1}\right)^{\frac{\lambda }{\widetilde{C}}}\left\vert f(x)\right\vert \left(
1+x\right) ^{\alpha }dx \\
&=&\frac{1}{\left( 1+y\right) ^{\alpha }r(y)}\int_{0}^{y}\left(\frac{x+1}{y+1}\right)^{%
\frac{\lambda }{\widetilde{C}}-\alpha }\left\vert f(x)\right\vert \left(
1+x\right) ^{\alpha }dx \\
&&\hspace{-2cm}\leq \frac{1}{\left( 1+y\right) ^{\alpha }r(y)}\int_{0}^{y}\left\vert
f(x)\right\vert \left( 1+x\right) ^{\alpha }dx\leq \frac{1}{\left(
1+y\right) ^{\alpha }r(y)}\left\Vert f\right\Vert _{X_{0,\alpha }},
\end{eqnarray*}%
on account of $\frac{x+1}{y+1}\leq 1$ and $\frac{\lambda }{\widetilde{C}}%
-\alpha \geq 0.$
\end{proof}

\subsection{The first perturbed semigroup}\mbox{}

We solve
\begin{equation*}
\ \frac{\partial }{\partial t}u(x,t)+\frac{\partial }{\partial x}\left[
r(x)u(x,t)\right] +a(x)u(x,t)=0
\end{equation*}%
by the method of characteristics. The solution is given by%
\begin{equation*}
\ e^{-\int_{X(y,t)}^{y}\frac{a(p)}{r(p)}dp}\frac{r(X(y,t))f(X(y,t))}{r(y)}.
\end{equation*}%
This defines a perturbed $C_{0}$-semigroup $\left( U(t)\right) _{t\geq 0}$
on both $X_{\alpha }$ and $X_{0,\alpha },\ $dominated by $\left(
U_{0}(t)\right) _{t\geq 0},$
\begin{equation*}
U(t)f=\ e^{-\int_{X(y,t)}^{y}\frac{a(p)}{r(p)}dp}\frac{r(X(y,t))f(X(y,t))}{%
r(y)}=e^{-\int_{X(y,t)}^{y}\frac{a(p)}{r(p)}dp}U_{0}(t)f.
\end{equation*}%
As previously, the Laplace transform of $(U(t))_{t\geqslant 0}$ and some
change of variables give:

\begin{proposition}
Let $\alpha >0$, (\ref{Main assumption Bis}) and (\ref{Assumptions on the
function r}) (resp. (\ref{Assumption on r bis})) be satisfied. The resolvent
of the generator $T$ of $\left( U(t)\right) _{t\geq 0}$ in $X_{\alpha }$
(resp. in $X_{0,\alpha }$), $\lambda >s(T)$, is given by
\begin{equation*}
\left( \left( \lambda -T\right) ^{-1}f\right) (y)=\frac{1}{r(y)}%
\int_{0}^{y}e^{-\int_{x}^{y}\frac{\lambda +\beta (p)}{r(\tau )}d\tau }f(x)dx.
\end{equation*}
\end{proposition}

\subsection{A smoothing effect of the perturbed resolvent}\mbox{}

The second a priori estimate in the spaces$\ X_{\alpha }$ and $X_{0,\alpha }$
is given by:

\begin{lemma}
\label{Lemma effet regularisant Bis}Let $\alpha >0$\ and (\ref{Main
assumption Bis}) be satisfied.

(i) Let (\ref{Assumptions on the function r}) be satisfied and $\lambda \geq
\alpha C.\ $Then, for any $f\in X_{\alpha },$
\begin{equation}
\int_{0}^{+\infty }\left\vert \left( (\lambda -T)^{-1}f\right)
(y)\right\vert a(y)y^{\alpha }dy\leq \ \int_{0}^{+\infty }\left\vert
(f(y)\right\vert y^{\alpha }dy.  \label{Smoothing effect 1}
\end{equation}

(ii) Let (\ref{Stronger assumption}) be satisfied and$\ \lambda \geq \alpha
\widetilde{C}.\ $Then, for any $f\in X_{0,\alpha },$
\begin{equation}
\int_{0}^{+\infty }\left\vert \left( (\lambda -T)^{-1}f\right)
(y)\right\vert a(y)\left( 1+y\right) ^{\alpha }dy\leq \ \int_{0}^{+\infty
}\left\vert (f(y)\right\vert \left( 1+y\right) ^{\alpha }dy.
\label{Smoothing effect 2}
\end{equation}
\end{lemma}

\begin{proof}
(i) By using (\ref{Estimate 2}) and $f\geq 0$
\begin{eqnarray*}
&&\int_{0}^{+\infty }\left( (\lambda -T)^{-1}f\right) (y)a(y)y^{\alpha }dy \\
&=&\int_{0}^{+\infty }\frac{a(y)y^{\alpha }}{r(y)}\left(
\int_{0}^{y}e^{-\lambda \int_{x}^{y}\frac{1}{r(p)}dp}e^{-\int_{x}^{y}\frac{%
a(p)}{r(p)}dp}f(x)dx\right) dy \\
&\leq &\int_{0}^{+\infty }\frac{a(y)y^{\alpha }}{r(y)}\left( \int_{0}^{y}
\left(\frac{x}{y}\right)^{\frac{\lambda }{C}}e^{-\int_{x}^{y}\frac{a(p)}{r(p)}%
dp}f(x)dx\right) dy \\
&=&\int_{0}^{+\infty }\left[ \int_{x}^{+\infty }\left(\frac{x}{y}\right)^{\frac{\lambda
}{C}}\frac{a(y)y^{\alpha }}{r(y)}e^{-\int_{x}^{y}\frac{a(p)}{r(p)}dp}dy%
\right] f(x)dx \\
&=&\int_{0}^{+\infty }\left[ \int_{x}^{+\infty }\frac{1}{x^{\alpha }}\left(\frac{x%
}{y}\right)^{\frac{\lambda }{C}}\frac{a(y)y^{\alpha }}{r(y)}e^{-\int_{x}^{y}\frac{%
a(p)}{r(p)}dp}dy\right] f(x)x^{\alpha }dx \\
&=&\int_{0}^{+\infty }\left[ \int_{x}^{+\infty }\left(\frac{x}{y}\right)^{\frac{\lambda
}{C}-\alpha }\frac{a(y)}{r(y)}e^{-\int_{x}^{y}\frac{a(p)}{r(p)}dp}dy\right]
f(x)x^{\alpha }dx \\
&\leq &\int_{0}^{+\infty }\left[ \int_{x}^{+\infty }\frac{a(y)}{r(y)}%
e^{-\int_{x}^{y}\frac{a(p)}{r(p)}dp}dy\right] f(x)x^{\alpha }dx,
\end{eqnarray*}%
where $\frac{x}{y}\leq 1$ and $\frac{\lambda }{C}-\alpha \geq 0\ $are used
in the last step. Thus%
\begin{eqnarray*}
&&\int_{0}^{+\infty }\left( (\lambda -T)^{-1}f\right) (y)a(y)y^{\alpha }dy \\
&\leq &\sup_{x>0}\int_{x}^{+\infty }\frac{a(y)}{r(y)}e^{-\int_{x}^{y}\frac{%
a(p)}{r(p)}dp}dy\left( \int_{0}^{+\infty }f(x)x^{\alpha }dx\right) .
\end{eqnarray*}%
Finally,
\begin{align*}
\int_{x}^{+\infty }e^{-\int_{x}^{y}\frac{a(p)}{r(p)}dp}\frac{a(y)}{r(y)}%
dy&=-\int_{x}^{+\infty }\frac{d}{dy}\left( e^{-\int_{x}^{y}\frac{a(p)}{r(p)}%
dp}\right) dy\\&=-\left[ e^{-\int_{x}^{y}\frac{a(p)}{r(p)}dp}\right]
_{y=x}^{y=+\infty }\leq 1
\end{align*}%
ends the proof.

(ii) The proof of (\ref{Smoothing effect 2}) is the same as that of Lemma %
\ref{Lemma effet regularisant} by using (\ref{Estimate 2 Bis}). \
\end{proof}

\subsection{On the full semigroup}\mbox{}

\textbf{\ }The same proof as that of Theorem \ref{Thm generation full
semigroup} in $X_{0,\alpha }$ gives the following statement which,
unfortunately, is \textit{not} useful for the purpose of spectral gaps, see
Remark \ref{Remark Open problem} below.

\begin{theorem}
Let (\ref{Main assumption Bis}) and (\ref{Assumption on r bis}) be
satisfied. Define%
\begin{equation*}
n_{1,\alpha }(y):=\int_{0}^{y}\left( 1+x\right) ^{\alpha }b(x,y)dx.
\end{equation*}%
If $\sup_{y>0}\frac{n_{1,\alpha }(y)}{\left( 1+y\right) ^{\alpha }}<+\infty
, $ then $B$ is $T$-bounded in $X_{0,\alpha }$ and
\begin{equation*}
\lim_{\lambda \rightarrow +\infty }\left\Vert B(\lambda -T)^{-1}\right\Vert
_{\mathcal{L}(X_{0,\alpha })}\leq \lim \sup_{a(y)\rightarrow +\infty }\frac{%
n_{1,\alpha }(y)}{\left( 1+y\right) ^{\alpha }}.
\end{equation*}%
In particular, if
\begin{equation}
\limsup_{a(y)\rightarrow +\infty }\frac{n_{1,\alpha }(y)}{\left(
1+y\right) ^{\alpha }}<1,  \label{Limsup less than 1}
\end{equation}%
then\textbf{\ }%
\begin{equation*}
A:=T+B: X_{0,\alpha }\supset D(T)\rightarrow X_{0,\alpha }
\end{equation*}%
generates a positive $C_{0}$-semigroup\textbf{\ }$\left( V(t)\right) _{t\geq
0}\ $on $X_{0,\alpha }.$
\end{theorem}

\begin{remark}
\label{Remark Open problem}If $a(\cdot)$ is\textit{\ unbounded near zero}, then
\begin{equation*}
\limsup_{a(y)\rightarrow +\infty }\frac{n_{1,\alpha }(y)}{\left(
1+y\right) ^{\alpha }}\geq \limsup_{y\rightarrow 0}n_{1,\alpha }(y)\geq
\limsup_{y\rightarrow 0}\int_{0}^{y}b(x,y)dx\geq 1,
\end{equation*}%
because%
\begin{equation*}
\int_{0}^{y}b(x,y)dx=\frac{1}{y}\int_{0}^{y}yb(x,y)dx\geq \frac{1}{y}%
\int_{0}^{y}xb(x,y)dx=1
\end{equation*}%
and hence (\ref{Limsup less than 1}) cannot be satisfied. On the
other hand, the compactness result we need in the sequel demands the
unboundedness of $a(\cdot)$ near zero. Hence, under Assumption (\ref{Main
assumption Bis}), it is not possible to finalize our spectral gap construction in
the space $X_{0,\alpha }.$ We point out that even if $a(\cdot)$ is\textit{\
unbounded near zero, we can still define a } positive $C_{0}$-semigroup%
\textbf{\ }$\left( V(t)\right) _{t\geq 0}\ $on $X_{0,\alpha }$ which solve
the growth fragmentation equations but in some generalized sense (honesty
theory), where the domain of the generator $T_{B}$ is the closure of $T+B$
only, see \cite[Chapter 5]{BLL}.\ However, in this case, we cannot infer that
$T_{B}$ is resolvent compact when $T$ is and the key argument behind the
existence of the spectral gap fails.
\end{remark}

\begin{theorem}
\label{Thm generation in X alpha}Let $\alpha >0$, (\ref{Main assumption Bis}%
) and (\ref{Assumptions on the function r}) be satisfied. Define%
\begin{equation*}
n_{\alpha }(y):=\int_{0}^{y}x^{\alpha }b(x,y)dx.
\end{equation*}%
If $\sup_{y>0}\frac{n_{\alpha }(y)}{y^{\alpha }}<+\infty ,$ then $B$ is $T$%
-bounded in $X_{\alpha }$ and
\begin{equation*}
\lim_{\lambda \rightarrow +\infty }\left\Vert B(\lambda -T)^{-1}\right\Vert
_{\mathcal{L}(X_{\alpha })}\leq \limsup_{a(y)\rightarrow +\infty }\frac{%
n_{\alpha }(y)}{y^{\alpha }}.
\end{equation*}%
In particular, if
\begin{equation}
\limsup_{a(y)\rightarrow +\infty }\frac{n_{\alpha }(y)}{y^{\alpha }}<1,
\label{Limsup less than 1 bis}
\end{equation}%
then\textbf{\ }%
\begin{equation*}
A:=T+B:D(T)\subset X_{\alpha }\rightarrow X_{\alpha }
\end{equation*}%
generates a positive $C_{0}$-semigroup\textbf{\ }$\left( V(t)\right) _{t\geq
0}\ $on $X_{\alpha }.$
\end{theorem}

\begin{proof}
We note that for nonnegative $\varphi, $ standard calculations give
\begin{equation*}
\left\Vert B\varphi \right\Vert _{X_{\alpha }}=\int_{0}^{+\infty
}a(y)n_{\alpha }(y)\varphi (y)dy.
\end{equation*}%
$\ $Thus, for nonnegative $f,$%
\begin{eqnarray}
\left\Vert B(\lambda -T)^{-1}f\right\Vert _{X_{\alpha }}  &=&\int_{0}^{+\infty }a(y)n_{\alpha }(y)\left( (\lambda -T)^{-1}f\right)
(y)dy  \notag \\
&=&\int_{0}^{+\infty }a(y)\frac{n_{\alpha }(y)}{y^{\alpha }}\left( (\lambda
-T)^{-1}f\right) (y)y^{\alpha }dy.  \label{To be decomposed bis}
\end{eqnarray}%
Let%
\begin{equation*}
L:=\limsup_{a(y)\rightarrow +\infty }\frac{n_{\alpha }(y)}{y^{\alpha }}.
\end{equation*}%
For any $\varepsilon >0$ there exists $c_{\varepsilon }>0$ such that%
\begin{equation*}
a(y)\geq c_{\varepsilon }\Longrightarrow \frac{n_{\alpha }(y)}{y^{\alpha }}%
\leq L+\varepsilon .
\end{equation*}%
We split (\ref{To be decomposed bis}) into two integrals
\begin{eqnarray*}
&&\int_{0}^{+\infty }a(y)\frac{n_{\alpha }(y)}{y^{\alpha }}\left( (\lambda
-T)^{-1}f\right) (y)y^{\alpha }dy \\
&=&\int_{\left\{ a(y)\leq c_{\varepsilon }\right\} }a(y)\frac{n_{\alpha }(y)%
}{y^{\alpha }}\left( (\lambda -T)^{-1}f\right) (y)y^{\alpha }dy \\
&&+\int_{\left\{ a(y)>c_{\varepsilon }\right\} }a(y)\frac{n_{\alpha }(y)}{%
y^{\alpha }}\left( (\lambda -T)^{-1}f\right) (y)y^{\alpha }dy \\
&=&I_{1}+I_{2}.
\end{eqnarray*}%
We note that
\begin{equation*}
I_{1}\leq c_{\varepsilon }\left\Vert \frac{n_{\alpha }(\cdot)}{y^{\alpha }}%
\right\Vert _{L^{\infty }}\left\Vert (\lambda -T)^{-1}f\right\Vert
_{X_{\alpha }}
\end{equation*}%
while, using Lemma \ref{Lemma effet regularisant Bis},%
\begin{eqnarray*}
I_{2} &\leq &\left( L+\varepsilon \right) \int_{0}^{+\infty }a(y)\left(
(\lambda -T)^{-1}f\right) (y)y^{\alpha }dy \\
&\leq &\left( L+\varepsilon \right) \left\Vert f\right\Vert _{X_{\alpha }}.
\end{eqnarray*}%
Hence,
\begin{eqnarray*}
\left\Vert B(\lambda -T)^{-1}f\right\Vert _{X_{\alpha }} &\leq
&c_{\varepsilon }\left\Vert \frac{n_{\alpha }(\cdot)}{y^{\alpha }}\right\Vert
_{L^{\infty }}\left\Vert (\lambda -T)^{-1}\right\Vert _{\mathcal{L}%
(X_{\alpha })}\left\Vert f\right\Vert _{X_{\alpha }} \\
&&+\left( L+\varepsilon \right) \left\Vert f\right\Vert _{X_{\alpha }}
\end{eqnarray*}%
and%
\begin{equation*}
\left\Vert B(\lambda -T)^{-1}\right\Vert _{\mathcal{L}(X_{\alpha })}\leq
c_{\varepsilon }\left\Vert \frac{n_{\alpha }(\cdot)}{y^{\alpha }}\right\Vert
_{L^{\infty }}\left\Vert (\lambda -T)^{-1}\right\Vert _{\mathcal{L}%
(X_{\alpha })}+\left( L+\varepsilon \right) .
\end{equation*}%
Since $\left\Vert (\lambda -T)^{-1}\right\Vert _{\mathcal{L}(X_{\alpha
})}\rightarrow 0\ $as $\lambda \rightarrow +\infty $,
\begin{equation*}
\lim_{\lambda \rightarrow +\infty }\left\Vert B(\lambda -T)^{-1}\right\Vert
_{\mathcal{L}(X_{\alpha })}\leq L+\varepsilon \ \ (\forall \varepsilon >0)
\end{equation*}%
and consequently, if $L<1,$ then
{\begin{equation*}
\lim_{\lambda \rightarrow +\infty }\left\Vert B(\lambda -T)^{-1}\right\Vert
_{\mathcal{L}(X_{\alpha })}< 1
\end{equation*}%
and we end the proof by applying Theorem \ref{Theorem Desch}}. $\ $
\end{proof}

\begin{remark}
\label{Remark monotony}{Proposition \ref{Prop monotony in alpha} yields that for each $y>0$, $\alpha \to \frac{n_{\alpha}(y)}{y^{\alpha}}$ is decreasing and convex. Since
\begin{equation*}
\frac{n_{1}(y)}{y}=\frac{1}{y}\int_{0}^{y}xb(x,y)dx=1,
\end{equation*}%
 $\limsup_{a(y)\rightarrow +\infty }\frac{n_{1}(y)}{y}=1 $ and
\begin{equation*}
\limsup\limits_{a(y)\rightarrow +\infty }\frac{n_{\alpha }(y)}{y^{\alpha }}\geq
1\ \ (0<\alpha \leq 1),
\end{equation*}%
hence the necessity to consider higher moments, i.e.,$\ \alpha >1$. So let $1<\alpha\leq \alpha_2$. By the convexity,
\begin{align*}
\frac{n_{\alpha}(y)}{y^\alpha} &\leq \frac{n_{1}(y)}{y} + \frac{\frac{n_{\alpha_2}(y)}{y^{\alpha_2}} -\frac{n_{1}(y)}{y}}{\alpha_2-1}(\alpha-1)= 1 + \frac{\frac{n_{\alpha_2}(y)}{y^{\alpha_2}} -1}{\alpha_2-1}(\alpha-1).
\end{align*}
If \eqref{Limsup less than 1 bis} holds for $\alpha_2$, for any $\epsilon >0,$ there is $c_\epsilon$ such that $\frac{n_{\alpha_2}(y)}{y^{\alpha_2}}\leq 1-\epsilon$ on   $\{y\in (0,\infty);\; a(y)\geq c_\epsilon\}$ and hence on this set
$$
\frac{n_{\alpha}(y)}{y^\alpha} \leq 1 -\epsilon\frac{\alpha-1}{\alpha_2-1}.
$$
Thus it follows that if (\ref{Limsup less than 1 bis}) holds for some $\alpha_2 >1,$
then it holds for all $\alpha>1.$}
\end{remark}

{
As noted in Introduction, for \textit{homogeneous} fragmentation kernels we have
\begin{equation*}
\frac{\int_{0}^{y}x^{\alpha }b(x,y)dx}{y^{\alpha }}=\int_{0}^{1}z^{\alpha
}h(z)dz<\int_{0}^{1}zh(z)dz=1\ (\alpha >1)
\end{equation*}%
\textbf{\ }and so\textbf{\ }$\limsup_{a(y)\rightarrow +\infty }\frac{n_{\alpha
}(y)}{y^{\alpha }}=\int_{0}^{1}z^{\alpha }h(z)dz<1$ for all $\alpha >1.$}

See also Section \ref{Section Separable models} for more examples.

\subsection{Compactness results}\mbox{}

By using Lemma \ref{Lemma pointwise estimate Bis}, Lemma \ref{Lemma effet
regularisant Bis} and arguing as in the proof of Theorem \ref{Thm T
resolvent compact} we get:

\begin{theorem}
Let $\alpha >0,\ $(\ref{Main assumption Bis}) and (\ref{Assumptions on the
function r})\ (resp. (\ref{Stronger assumption})) be satisfied. Let the
sublevel sets of $a(\cdot)$ be thin at zero and at infinity in the sense that
for any $c>0$
\begin{equation*}
\int_{0}^{+\infty }1_{\left\{ a<c\right\} }\frac{1\ }{r(y)}dy<+\infty \
\end{equation*}%
(e.g., let $\lim_{x\rightarrow +\infty }a(x)=+\infty $ and $%
\lim_{x\rightarrow 0}a(x)=+\infty $). Then $T$ is resolvent compact in $%
X_{\alpha }$ (resp. in $X_{0,\alpha }$).
\end{theorem}

Similarly to Corollary \ref{Corollary resolvent compact} we have
\begin{corollary}
Let (\ref{Main assumption Bis}), (\ref{Assumptions on the function r}) and (%
\ref{Limsup less than 1 bis}) be satisfied.\ If the sublevel sets of $a(\cdot)$
are thin at zero and at infinity (e.g. if $\ \lim_{x\rightarrow +\infty
}a(x)=+\infty $ and $\lim_{x\rightarrow 0}a(x)=+\infty $), then $A:=T+B\ $ is
resolvent compact in $X_{\alpha }$.\
\end{corollary}

\subsection{Spectral gap of $(V(t))_{t\geqslant 0}$ in $X_{\protect\alpha }$}\mbox{}

The same proof as for Theorem \ref{irred} gives:

\begin{lemma}
{Let $\alpha>0$ and assumptions 1. or 2. of Theorem \ref{irred}},\ (\ref{Main assumption Bis}), (\ref%
{Assumptions on the function r}) and (\ref{Limsup less than 1 bis}) be
satisfied. Then $(\lambda -T-B)^{-1}$ is positivity improving in $X_{\alpha
}\ $or, equivalently, the $C_{0}$-semigroup $\left( V(t)\right) _{t\geq 0}$
is irreducible in $X_{\alpha }$.
\end{lemma}

Finally, the same proof as for Theorem \ref{Thm Main result} gives the main
result of the second construction.

\begin{theorem}
\label{Thm Main result Bis}Let $\alpha >0$ {and assumptions 1. or 2. of Theorem \ref{irred}},\ (\ref{Main assumption Bis}), (\ref{Assumptions on the function r}%
) and (\ref{Limsup less than 1 bis}) be satisfied. If the sublevel sets of $%
a(\cdot)$ are thin at zero and at infinity, (e.g. if $\ \lim_{x\rightarrow
+\infty }a(x)=+\infty $ and $\lim_{x\rightarrow 0}a(x)=+\infty $), then $%
\left( V(t)\right) _{t\geq 0}$\ has an asynchronous exponential growth in $%
X_{\alpha }.$
\end{theorem}

\begin{remark}
We conjecture that the result does not hold without the unboundedness of $%
a(\cdot)$ at zero as suggested by \cite[Theorem 4.1]{BG2}.
\end{remark}

\section*{Acknowledgments} This work began while the first author was visiting the University of
Pretoria\ in January-February 2020. The support for this visit and research came from the DSI/NRF SARChI Grant 82770. The second author was also supported by the National Science Centre of Poland Grant 2017/25/B/ST1/00051.

\appendix

\section{\label{Section Separable models}Separable fragmentation kernels}

We have seen how homogeneous fragmentation kernels satisfy the key
assumptions (\ref{Limit less than 1}) and (\ref{Limsup less than 1 bis}) of
our construction. This last section is devoted to separable fragmentation
kernels%
\begin{equation*}
b(x,y)=\beta (x)\gamma (y),
\end{equation*}%
introduced in \cite{Banasiak 2004}, see also \cite{Banasiak-Arlotti}\cite%
{BLL}. It is easy to see that separable kernels with mass conservation (\ref%
{conservativity of mass}) are of the form
\begin{equation}
b(x,y)=\beta (x)y\left( \int_{0}^{y}s\beta (s)ds\right) ^{-1}
\label{Separable kernel}
\end{equation}%
where
\begin{equation*}
0<\int_{0}^{y}s\beta (s)ds<+\infty \ \ \forall y>0.
\end{equation*}%
A particular case of separable kernels are power law kernels%
\begin{equation*}
b(x,y)=\left( \nu +2\right) \frac{x^{\nu }}{y^{\nu +1}}\ \ (\nu \in \left(
-2,0\right] ).
\end{equation*}
We can complement Theorem \ref{Thm generation in X alpha} by:

\begin{proposition}
We assume that the fragmentation kernel is of the form (\ref{Separable
kernel}) and $a(\cdot)$ is only unbounded at zero and infinity. If%
\begin{equation*}
\beta _{0}^{-}:=\lim \inf_{x\rightarrow 0}x\beta (x)>0\text{,\ }\beta
_{0}^{+}:=\limsup_{x\rightarrow 0}x\beta (x)<+\infty
\end{equation*}%
and%
\begin{equation*}
\int_{0}^{+\infty }x^{\alpha }\beta (x)dx<+\infty ,
\end{equation*}%
then (\ref{Limsup less than 1 bis}) is satisfied provided $\alpha >\frac{%
\beta _{0}^{+}}{\beta _{0}^{-}}.$
\end{proposition}

\begin{proof}
Note first that
\begin{equation*}
n_{\alpha }(y)=\left( \int_{0}^{y}x^{\alpha }\beta (x)dx\right) y\left(
\int_{0}^{y}x\beta (x)dx\right) ^{-1}
\end{equation*}%
and (\ref{Limsup less than 1 bis}) amounts to%
\begin{equation}
\limsup_{y\rightarrow +\infty }\frac{1}{y^{\alpha -1}}\frac{%
\int_{0}^{y}x^{\alpha }\beta (x)dx}{\int_{0}^{y}x\beta (x)dx}<1
\label{Limsup at infinity}
\end{equation}%
and%
\begin{equation}
\limsup_{y\rightarrow 0}\frac{1}{y^{\alpha -1}}\frac{\int_{0}^{y}x^{\alpha
}\beta (x)dx}{\int_{0}^{y}x\beta (x)dx}<1.  \label{Limsup at zero}
\end{equation}%
Let $\varepsilon >0$ be arbitrary.\ Then, for $y$ small enough%
\begin{equation*}
\int_{0}^{y}x^{\alpha }\beta (x)dx\leq \left( \beta _{0}^{+}+\varepsilon
\right) \int_{0}^{y}x^{\alpha -1}dx=\frac{\left( \beta _{0}^{+}+\varepsilon
\right) y^{\alpha }}{\alpha }
\end{equation*}%
and
\begin{equation*}
\int_{0}^{y}x\beta (x)dx\geq \left( \beta _{0}^{-}-\varepsilon \right) y,
\end{equation*}%
so%
\begin{equation*}
\frac{1}{y^{\alpha -1}}\frac{\int_{0}^{y}x^{\alpha }\beta (x)dx}{%
\int_{0}^{y}x\beta (x)dx}\leq \frac{1}{\alpha }\frac{\beta
_{0}^{+}+\varepsilon }{\beta _{0}^{-}-\varepsilon }.
\end{equation*}%
Therefore,
\begin{equation*}
\limsup_{y\rightarrow 0}\frac{1}{y^{\alpha -1}}\frac{\int_{0}^{y}x^{\alpha
}\beta (x)dx}{\int_{0}^{y}x\beta (x)dx}\leq \frac{1}{\alpha }\frac{\beta
_{0}^{+}}{\beta _{0}^{-}}.
\end{equation*}%
Finally,%
\begin{equation*}
\frac{1}{y^{\alpha -1}}\frac{\int_{0}^{y}x^{\alpha }\beta (x)dx}{%
\int_{0}^{y}x\beta (x)dx}\leq \frac{1}{y^{\alpha -1}}\frac{\int_{0}^{\infty
}x^{\alpha }\beta (x)dx}{\int_{0}^{1}x\beta (x)dx}\ (y\geq 1)
\end{equation*}%
and%
\begin{equation*}
\lim_{y\rightarrow +\infty }\frac{1}{y^{\alpha -1}}\frac{\int_{0}^{y}x^{%
\alpha }\beta (x)dx}{\int_{0}^{y}x\beta (x)dx}=0.
\end{equation*}%
This ends the proof.
\end{proof}

\begin{remark}
If $\beta _{0}:=\lim_{x\rightarrow 0}x\beta (x)>0$ exists then $\beta
_{0}^{+}=\beta _{0}^{-}$ and both (\ref{Limsup at infinity}) and (\ref{Limsup at
zero}) are satisfied for any $\alpha >1$ such \ that $\int_{0}^{+\infty
}x^{\alpha }\beta (x)dx<+\infty .$
\end{remark}

Similarly, we can complement Theorem \ref{Thm generation full semigroup} by:

\begin{proposition}
We assume that the fragmentation kernel is of the form (\ref{Separable
kernel}) and $a(\cdot)$ is only unbounded at infinity. If\
\begin{equation}
\int_{0}^{+\infty }\beta (x)\left( 1+x^{\overline{\alpha }}\right) dx<+\infty
\label{Finite integral}
\end{equation}%
for some $\overline{\alpha }>1$ then (\ref{Limit less than 1}) is satisfied
for any $\alpha >1$ and consequently the \textit{threshold is equal to one.}
\end{proposition}

Note first that
\begin{equation*}
n_{1,\alpha }(y):=y\left( \int_{0}^{y}x\beta (x)dx\right)
^{-1}\int_{0}^{y}\beta (x)\left( 1+x\right) ^{\alpha }dx
\end{equation*}%
and hence (\ref{Limit less than 1}) amounts to%
\begin{equation*}
\limsup_{y\rightarrow +\infty }\frac{y}{\left( 1+y\right) ^{\alpha }}\frac{%
\int_{0}^{y}\beta (x)\left( 1+x\right) ^{\alpha }dx}{\int_{0}^{y}x\beta (x)dx%
}<1.
\end{equation*}%
Note that (\ref{Finite integral}) implies that
\begin{equation*}
\int_{0}^{+\infty }\beta (x)\left( 1+x^{\alpha }\right) dx<+\infty \ \
(0<\alpha \leq \overline{\alpha }).
\end{equation*}%
It is easy to see that $\left( 1+x\right) ^{\alpha }\leq 2^{\alpha -1}\left(
1+x^{\alpha }\right) $ so%
\begin{equation*}
\int_{0}^{y}\beta (x)\left( 1+x\right) ^{\alpha }dx\leq 2^{\alpha
-1}\int_{0}^{y}\beta (x)\left( 1+x^{\alpha }\right) dx
\end{equation*}%
and consequently, for any $1<\alpha \leq \overline{\alpha },$
\begin{equation*}
\frac{y}{\left( 1+y\right) ^{\alpha }}\frac{\int_{0}^{y}\beta (x)\left(
1+x\right) ^{\alpha }dx}{\int_{0}^{y}x\beta (x)dx}\leq \frac{1}{y^{\alpha -1}%
}\frac{2^{\alpha -1}\int_{0}^{y}\beta (x)\left( 1+x^{\alpha }\right) dx}{%
\int_{0}^{y}x\beta (x)dx}\rightarrow 0\ \ (y\rightarrow +\infty )
\end{equation*}%
and this ends the proof.

\begin{remark}
We note that any convex combination of conservative fragmentation kernels is
a conservative fragmentation kernel so
\begin{equation}
b(x,y)=\sum_{j\in J}\lambda _{j}\beta _{j}(x)y\left( \int_{0}^{y}s\beta
_{j}(s)ds\right) ^{-1},\ \left(\sum_{j\in J}\lambda _{j}=1\right)
\label{Degenerate fragmentation kernel}
\end{equation}%
($J$ finite or denumerable) is a conservative kernel and we can check the
key conditions $\limsup_{a(y)\rightarrow +\infty }\frac{n_{\alpha }(y)}{%
y^{\alpha }}<1$ or $\limsup_{a(y)\rightarrow +\infty }\frac{n_{1,\alpha
}(y)}{\left( 1+y\right) ^{\alpha }}<1$ more generally for (\ref{Degenerate
fragmentation kernel}) by using just the last two propositions above. We
could also consider convex combinations of separable kernels and homogeneous
ones.
\end{remark}

\end{document}